# SUBTREE PRUNE AND REGRAFT: A REVERSIBLE REAL TREE-VALUED MARKOV PROCESS

By Steven N. Evans[1] and Anita Winter

*University of California at Berkeley and Universität Erlangen–Nürnberg*

We use Dirichlet form methods to construct and analyze a reversible Markov process, the stationary distribution of which is the Brownian continuum random tree. This process is inspired by the subtree prune and regraft (SPR) Markov chains that appear in phylogenetic analysis.

A key technical ingredient in this work is the use of a novel Gromov–Hausdorff type distance to metrize the space whose elements are compact real trees equipped with a probability measure. Also, the investigation of the Dirichlet form hinges on a new path decomposition of the Brownian excursion.

**1. Introduction.** Markov chains that move through a space of finite trees are an important ingredient for several algorithms in phylogenetic analysis, particularly in Markov chain Monte Carlo algorithms for simulating distributions on spaces of trees in Bayesian tree reconstruction and in simulated annealing algorithms in maximum likelihood and maximum parsimony tree reconstruction (see, e.g., [21] for a comprehensive overview of the field). (*Maximum parsimony* tree reconstruction is based on finding the phylogenetic tree and inferred ancestral states that minimize the total number of obligatory inferred substitution events on the edges of the tree.) Usually, such chains are based on a set of simple rearrangements that transform a tree into a "neighboring" tree. One widely used set of moves is the *nearest-neighbor interchanges* (NNI) (see, e.g., [6, 7, 9, 21]). Two other standard sets of moves that are implemented in several phylogenetic software packages but seem to have received less theoretical attention are the *subtree prune and regraft* (SPR) moves and the *tree bisection and reconnection* (TBR) moves

Received February 2005; revised August 2005.

[1]Supported in part by NSF Grants DMS-00-71468 and DMS-04-05778.

that were first described in [32] and are further discussed in [6, 21, 30]. We note that an NNI move is a particular type of SPR move and that an SPR move is a particular type of TBR move and, moreover, that every TBR operation is either a single SPR move or the composition of two such moves (see, e.g., Section 2.6 of [30]). Chains based on other moves are investigated in [5, 14, 29].

In an SPR move, a binary tree $T$ (i.e., a tree in which all nonleaf vertices have degree 3) is cut "in the middle of an edge" to give two subtrees, say $T'$ and $T''$. Another edge is chosen in $T'$, a new vertex is created "in the middle" of that edge and the cut edge in $T''$ is attached to this new vertex. Last, the "pendant" cut edge in $T'$ is removed along with the vertex it was attached to in order to produce a new binary tree that has the same number of vertices as $T$. See Figure 1.

As remarked in [6],

> The SPR operation is of particular interest as it can be used to model biological processes such as horizontal gene transfer and recombination.

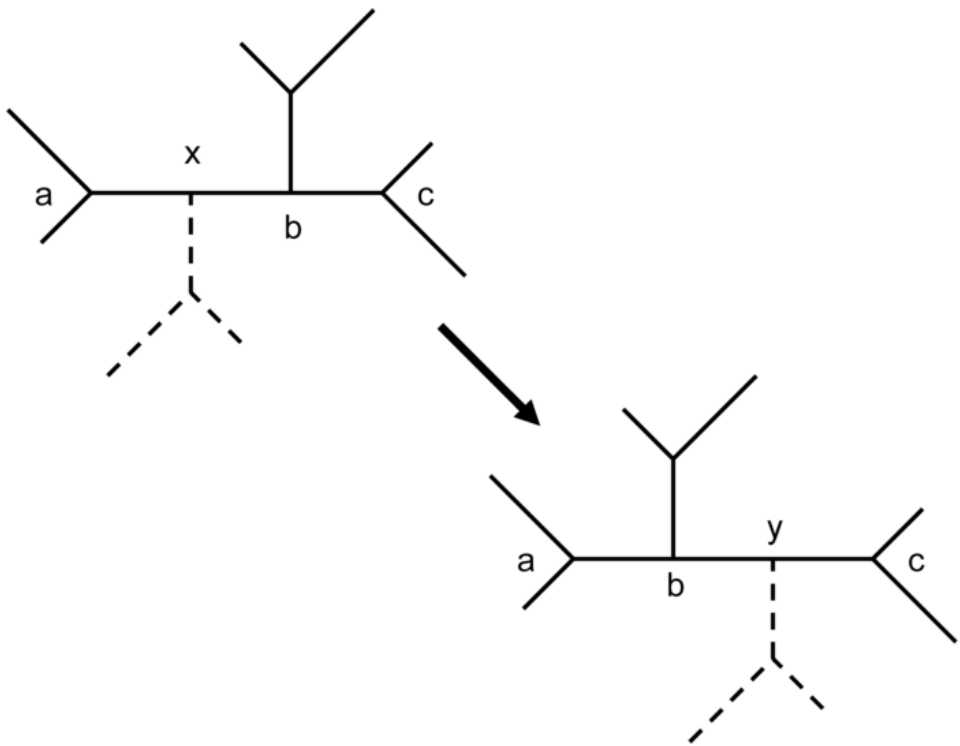

FIG. 1. *An SPR move. The dashed subtree tree attached to vertex $x$ in the top tree is reattached at a new vertex $y$ that is inserted into the edge $(b,c)$ in the bottom tree to make two edges $(b,y)$ and $(y,c)$. The two edges $(a,x)$ and $(b,x)$ in the top tree are merged into a single edge $(a,b)$ in the bottom tree.*

(*Horizontal gene transfer* is the transfer of genetic material from one species to another. It is a particularly common phenomenon among bacteria.) Section 2.7 of [30] provides more background on this point as well as a comment on the role of SPR moves in the two phenomena of lineage sorting and gene duplication and loss.

In this paper we investigate the asymptotics of the simplest possible tree-valued Markov chain based on the SPR moves, namely the chain in which the two edges that are chosen for cutting and for reattaching are chosen uniformly (without replacement) from the edges in the current tree. Intuitively, the continuous-time Markov process we discuss arises as limit when the number of vertices in the tree goes to infinity, the edge lengths are rescaled by a constant factor so that the initial tree converges in a suitable sense to a continuous analogue of a combinatorial tree (more specifically, a compact real tree), and the time scale of the Markov chain is sped up by an appropriate factor.

We do not, in fact, prove such a limit theorem. Rather, we use Dirichlet form techniques to establish the existence of a process that has the dynamics one would expect from such a limit. Unfortunately, although Dirichlet form techniques provide powerful tools for constructing and analyzing symmetric Markov processes, they are notoriously inadequate for proving convergence theorems (as opposed to generator or martingale problem characterizations of Markov processes, e.g.). We therefore leave the problem of establishing a limit theorem to future research.

The Markov process we construct is a pure jump process that is reversible with respect to the distribution of Aldous' continuum random tree (i.e., the random tree which arises as the rescaling limit of uniform random trees with $n$ vertices when $n \to \infty$ and which is also, up to a constant scaling factor, the random tree associated naturally with the standard Brownian excursion—see Section 4 for more details about the continuum random tree, its connection with Brownian excursion and references to the literature).

Somewhat more precisely, but still rather informally, the process we construct has the following description.

To begin with, Aldous' continuum random tree has two natural measures on it that can both be thought of as arising from the measure on an approximating finite tree with $n$ vertices that places a unit mass at each vertex. If we rescale the mass of this measure to get a probability measure, then in the limit we obtain a probability measure on the continuum random tree that happens to assign all of its mass to the leaves with probability 1. We call this probability measure the *weight* on the continuum tree. On the other hand, we can also rescale the measure that places a unit mass at each vertex to obtain in the limit a $\sigma$-finite measure on the continuum tree that restricts to one-dimensional Lebesgue measure if we restrict to any path through the continuum tree. We call this $\sigma$-finite measure the *length*.

The continuum random tree is a random compact *real tree* of the sort investigated in [20] (we define real trees and discuss some of their properties in Section 2). Any compact real tree has an analogue of the length measure on it, but in general there is no canonical analogue of the weight measure. Consequently, the process we construct has as its state space the set of pairs $(T, \nu)$, where $T$ is a compact real tree and $\nu$ is a probability measure on $T$. Let $\mu$ be the length measure associated with $T$. Our process jumps away from $T$ by first choosing a pair of points $(u, v) \in T \times T$ according to the rate measure $\mu \otimes \nu$ and then transforming $T$ into a new tree by cutting off the subtree rooted at $u$ that does not contain $v$ and reattaching this subtree at $v$. This jump kernel (which typically has infinite total mass—so that jumps are occurring on a dense countable set) is precisely what one would expect for a limit (as the number of vertices goes to infinity) of the particular SPR Markov chain on finite trees described above in which the edges for cutting and reattachment are chosen uniformly at each stage.

The framework of Dirichlet forms allows us to translate this description into rigorous mathematics. An important preliminary step that we accomplish in Section 2 is to show that it is possible to equip the space of pairs of compact real trees and their accompanying weights with a nice Gromov–Hausdorff-like metric that makes this space complete and separable. We note that a Gromov–Hausdorff-like metric on more general metric spaces equipped with measures was introduced in [31]. The latter metric is based on the Wasserstein $L^2$ distance between measures, whereas ours is based on the Prohorov distance. Moreover, we need to understand in detail the Dirichlet form arising from the combination of the jump kernel with the continuum random tree distribution as a reference measure, and we accomplish this in Sections 5 and 6, where we establish the relevant facts from what appears to be a novel path decomposition of the standard Brownian excursion. We construct the Dirichlet form and the resulting process in Section 7. We use potential theory for Dirichlet forms to show in Section 8 that from almost all starting points (with respect to the continuum random tree reference measure) our process does not hit the trivial tree consisting of a single point.

We remark that excursion path-valued Markov processes that are reversible with respect to the distribution of standard Brownian excursion and have *continuous* sample paths have been investigated in [34, 35, 36], and that these processes can also be thought of as real tree-valued diffusion processes that are reversible with respect to the distribution of the continuum random tree. However, we are unaware of a description in which these latter processes arise as limits of natural processes on spaces of finite trees.

**2. Weighted $\mathbb{R}$-trees.** A metric space $(X, d)$ is a *real tree* ($\mathbb{R}$-tree) if it satisfies the following axioms.

AXIOM 0 (Completeness). The space $(X, d)$ is complete.

AXIOM 1 (Unique geodesics). For all $x, y \in X$ there exists a unique isometric embedding $\phi_{x,y} \colon [0, d(x, y)] \to X$ such that $\phi_{x,y}(0) = x$ and $\phi_{x,y}(d(x, y)) = y$.

AXIOM 2 (Loop-free). For every injective continuous map $\psi \colon [0, 1] \to X$ one has $\psi([0, 1]) = \phi_{\psi(0), \psi(1)}([0, d(\psi(0), \psi(1))])$.

Axiom 1 says simply that there is a unique "unit speed" path between any two points, whereas Axiom 2 implies that the image of any injective path connecting two points coincides with the image of the unique unit speed path, so that it can be reparameterized to become the unit speed path. Thus, Axiom 1 is satisfied by many other spaces such as $\mathbb{R}^d$ with the usual metric, whereas Axiom 2 expresses the property of "treeness" and is only satisfied by $\mathbb{R}^d$ when $d = 1$. We refer the reader to [12, 15, 16, 17, 33] for background on $\mathbb{R}$-trees. In particular, [12] shows that a number of other definitions are equivalent to the one above. A particularly useful fact is that a metric space $(X, d)$ is an $\mathbb{R}$-tree if and only if it is complete, path-connected and satisfies the so-called *four point condition*, that is,

$$
\begin{aligned}
(2.1) \quad & d(x_1, x_2) + d(x_3, x_4) \\
& \leq \max\{d(x_1, x_3) + d(x_2, x_4), d(x_1, x_4) + d(x_2, x_3)\}
\end{aligned}
$$

for all $x_1, \ldots, x_4 \in X$.

Let $\mathbf{T}$ denote the set of isometry classes of compact $\mathbb{R}$-*trees*. In order to equip $\mathbf{T}$ with a metric, recall that the *Hausdorff distance* between two closed subsets $A$, $B$ of a metric space $(X, d)$ is defined as

$$(2.2) \qquad d_\mathrm{H}(A, B) := \inf\{\varepsilon > 0 \colon A \subseteq U_\varepsilon(B) \text{ and } B \subseteq U_\varepsilon(A)\},$$

where $U_\varepsilon(C) := \{x \in X \colon d(x, C) \leq \varepsilon\}$. Based on this notion of distance between closed sets, we define the *Gromov–Hausdorff distance*, $d_{\mathrm{GH}}(X, Y)$, between two metric spaces $(X, d_X)$ and $(Y, d_Y)$ as the infimum of the Hausdorff distance $d_\mathrm{H}(X', Y')$ over all metric spaces $X'$ and $Y'$ that are isomorphic to $X$ and $Y$, respectively, and that are subspaces of some common metric space $Z$ (cf. [10, 11, 23]).

A direct application of the previous definition requires an optimal embedding into a space $Z$ which it is not possible to obtain explicitly in most examples. We therefore give an equivalent reformulation which allows us to get estimates on the distance by looking for "matchings" between the two spaces that preserve the two metrics up to an additive error. In order to be more explicit, we require some more notation. A subset $\Re \subseteq X \times Y$ is said to be a *correspondence* between sets $X$ and $Y$ if for each $x \in X$ there exists at least one $y \in Y$ such that $(x,y) \in \Re$, and for each $y \in Y$ there exists at least one $x \in X$ such that $(x,y) \in \Re$. The *distortion* of $\Re$ is defined by

$$(2.3) \quad \operatorname{dis}(\Re) := \sup\{|d_X(x_1, x_2) - d_Y(y_1, y_2)| : (x_1, y_1), (x_2, y_2) \in \Re\}.$$

Then

$$(2.4) \qquad d_{\operatorname{GH}}((X, d_X), (Y, d_Y)) = \tfrac{1}{2} \inf_{\Re} \operatorname{dis}(\Re),$$

where the infimum is taken over all correspondences $\Re$ between $X$ and $Y$ (see, e.g., Theorem 7.3.25 in [11]).

It is shown in Theorem 1 in [20] that the metric space $(\mathbf{T}, d_{\operatorname{GH}})$ is complete and separable.

In the following we will be interested in compact $\mathbb{R}$-trees $(T,d) \in \mathbf{T}$ equipped with a probability measure $\nu$ on the Borel $\sigma$-field $\mathcal{B}(T)$. We call such objects weighted compact $\mathbb{R}$-trees and write $\mathbf{T}^{\operatorname{wt}}$ for the space of weight-preserving isometry classes of weighted compact $\mathbb{R}$-trees, where we say that two weighted, compact $\mathbb{R}$-trees $(X, d, \nu)$ and $(X', d', \nu')$ are weight-preserving isometric if there exists an isometry $\phi$ between $X$ and $X'$ such that the *push-forward* of $\nu$ by $\phi$ is $\nu'$:

$$(2.5) \qquad \nu' = \phi_* \nu := \nu \circ \phi^{-1}.$$

It is clear that the property of being weight-preserving isometric is an equivalence relation.

We want to equip $\mathbf{T}^{\operatorname{wt}}$ with a Gromov–Hausdorff type of distance which incorporates the weights on the trees, but first we need to introduce some notions that will be used in the definition.

An $\varepsilon$-(*distorted*) *isometry* between two metric spaces $(X, d_X)$ and $(Y, d_Y)$ is a (possibly nonmeasurable) map $f : X \to Y$ such that

$$(2.6) \quad \operatorname{dis}(f) := \sup\{|d_X(x_1, x_2) - d_Y(f(x_1), f(x_2))| : x_1, x_2 \in X\} \leq \varepsilon$$

and $f(X)$ is an $\varepsilon$-net in $Y$.

It is easy to see that if for two metric spaces $(X, d_X)$ and $(Y, d_Y)$ and $\varepsilon > 0$ we have $d_{\operatorname{GH}}((X, d_X), (Y, d_Y)) < \varepsilon$, then there exists a $2\varepsilon$-isometry from $X$ to $Y$ (cf. Lemma 7.3.28 in [11]). The following lemma states that we may choose the distorted isometry between $X$ and $Y$ to be measurable if we allow a slightly bigger distortion.

LEMMA 2.1.  *Let $(X, d_X)$ and $(Y, d_Y)$ be two compact real trees such that $d_{\operatorname{GH}}((X, d_X), (Y, d_Y)) < \varepsilon$ for some $\varepsilon > 0$. Then there exists a measurable $3\varepsilon$-isometry from $X$ to $Y$.*

PROOF.  If $d_{\operatorname{GH}}((X, d_X), (Y, d_Y)) < \varepsilon$, then by (2.4) there exists a correspondence $\Re$ between $X$ and $Y$ such that $\operatorname{dis}(\Re) < 2\varepsilon$. Since $(X, d_X)$ is compact there exists a finite $\varepsilon$-net in $X$. We claim that for each such finite $\varepsilon$-net $S^{X,\varepsilon} = \{x_1, \ldots, x_{N^\varepsilon}\} \subseteq X$, any set $S^{Y,\varepsilon} = \{y_1, \ldots, y_{N^\varepsilon}\} \subseteq Y$ such that $(x_i, y_i) \in \Re$ for all $i \in \{1, 2, \ldots, N^\varepsilon\}$ is a $3\varepsilon$-net in $Y$. To see this, fix $y \in Y$. We have to show the existence of $i \in \{1, 2, \ldots, N^\varepsilon\}$ with $d_Y(y_i, y) < 3\varepsilon$. For that choose $x \in X$ such that $(x, y) \in \Re$. Since $S^{X,\varepsilon}$ is an $\varepsilon$-net in $X$ there exists an $i \in \{1, 2, \ldots, N^\varepsilon\}$ such that $d_X(x_i, x) < \varepsilon$. $(x_i, y_i) \in \Re$ implies therefore that $|d_X(x_i, x) - d_Y(y_i, y)| \leq \text{dis}(\Re) < 2\varepsilon$, and hence $d_Y(y_i, y) < 3\varepsilon$.

Furthermore we may decompose $X$ into $N^\varepsilon$ possibly empty measurable disjoint subsets of $X$ by letting $X^{1,\varepsilon} := \mathcal{B}(x_1, \varepsilon)$, $X^{2,\varepsilon} := \mathcal{B}(x_2, \varepsilon) \setminus X^{1,\varepsilon}$, and so on, where $\mathcal{B}(x, r)$ is the open ball $\{x' \in X : d_X(x, x') < r\}$. Then $f$ defined by $f(x) = y_i$ for $x \in X^{i,\varepsilon}$ is obviously a measurable $3\varepsilon$-isometry from $X$ to $Y$. □

We also need to recall the definition of the Prohorov distance between two probability measures (see, e.g., [19]). Given two probability measures $\mu$ and $\nu$ on a metric space $(X, d)$ with the corresponding collection of closed sets denoted by $\mathcal{C}$, the Prohorov distance between them is

$$d_\text{P}(\mu, \nu) := \inf\{\varepsilon > 0 : \mu(C) \leq \nu(C^\varepsilon) + \varepsilon \text{ for all } C \in \mathcal{C}\},$$

where $C^\varepsilon := \{x \in X : \inf_{y \in C} d(x, y) < \varepsilon\}$. The Prohorov distance is a metric on the collection of probability measures on $X$. The following result shows that if we push measures forward with a map having a small distortion, then Prohorov distances cannot increase too much.

LEMMA 2.2. *Suppose that $(X, d_X)$ and $(Y, d_Y)$ are two metric spaces, $f : X \to Y$ is a measurable map with $\text{dis}(f) \leq \varepsilon$, and $\mu$ and $\nu$ are two probability measures on $X$. Then*

$$d_\text{P}(f_*\mu, f_*\nu) \leq d_\text{P}(\mu, \nu) + \varepsilon.$$

PROOF. Suppose that $d_\text{P}(\mu, \nu) < \delta$. By definition, $\mu(C) \leq \nu(C^\delta) + \delta$ for all closed sets $C \in \mathcal{C}$. If $D$ is a closed subset of $Y$, then

$$\begin{aligned} f_*\mu(D) &= \mu(f^{-1}(D)) \\ &\leq \mu(\overline{f^{-1}(D)}) \\ &\leq \nu(\overline{f^{-1}(D)}^\delta) + \delta \\ &= \nu(f^{-1}(D)^\delta) + \delta. \end{aligned} \quad (2.7)$$

Now $x' \in f^{-1}(D)^\delta$ means there is $x'' \in X$ such that $d_X(x', x'') < \delta$ and $f(x'') \in D$. By the assumption that $\text{dis}(f) \leq \varepsilon$, we have $d_Y(f(x'), f(x'')) < \delta + \varepsilon$, and hence $f(x') \in D^{\delta+\varepsilon}$. Thus

$$f^{-1}(D)^\delta \subseteq f^{-1}(D^{\delta+\varepsilon}) \quad (2.8)$$

and we have

$$f_*\mu(D) \leq \nu(f^{-1}(D^{\delta+\varepsilon})) + \delta = f_*\nu(D^{\delta+\varepsilon}) + \delta, \tag{2.9}$$

so that $d_{\mathrm{P}}(f_*\mu, f_*\nu) \leq \delta + \varepsilon$, as required. $\square$

We are now in a position to define the weighted Gromov–Hausdorff distance between the two compact, weighted $\mathbb{R}$-trees $(X, d_X, \nu_X)$ and $(Y, d_Y, \nu_Y)$. For $\varepsilon > 0$, set

$$F_{X,Y}^{\varepsilon} := \{\text{measurable } \varepsilon\text{-isometries from } X \text{ to } Y\}. \tag{2.10}$$

Put

$$\begin{aligned}
\Delta_{\mathrm{GH^{wt}}}(X, Y) \\
:= \inf\{\varepsilon > 0 : \text{exist } f \in F_{X,Y}^{\varepsilon}, g \in F_{Y,X}^{\varepsilon} \text{ such that} \\
d_{\mathrm{P}}(f_*\nu_X, \nu_Y) \leq \varepsilon, d_{\mathrm{P}}(\nu_X, g_*\nu_Y) \leq \varepsilon\}.
\end{aligned} \tag{2.11}$$

Note that the set on the right-hand side is nonempty because $X$ and $Y$ are compact, and hence bounded. It will turn out that $\Delta_{\mathrm{GH^{wt}}}$ satisfies all the properties of a metric except the triangle inequality. To rectify this, let

$$d_{\mathrm{GH^{wt}}}(X, Y) := \inf\left\{\sum_{i=1}^{n-1} \Delta_{\mathrm{GH^{wt}}}(Z_i, Z_{i+1})^{1/4}\right\}, \tag{2.12}$$

where the infimum is taken over all finite sequences of compact, weighted $\mathbb{R}$-trees $Z_1, \ldots, Z_n$ with $Z_1 = X$ and $Z_n = Y$.

LEMMA 2.3. *The map $d_{\mathrm{GH^{wt}}} : \mathbf{T}^{\mathrm{wt}} \times \mathbf{T}^{\mathrm{wt}} \to \mathbb{R}_+$ is a metric on $\mathbf{T}^{\mathrm{wt}}$. Moreover,*

$$\tfrac{1}{2}\Delta_{\mathrm{GH^{wt}}}(X, Y)^{1/4} \leq d_{\mathrm{GH^{wt}}}(X, Y) \leq \Delta_{\mathrm{GH^{wt}}}(X, Y)^{1/4}$$

*for all $X, Y \in \mathbf{T}^{\mathrm{wt}}$.*

PROOF. It is immediate from (2.11) that the map $\Delta_{\mathrm{GH^{wt}}}$ is symmetric. We next claim that

$$\Delta_{\mathrm{GH^{wt}}}((X, d_X, \nu_X), (Y, d_Y, \nu_Y)) = 0, \tag{2.13}$$

if and only if $(X, d_X, \nu_X)$ and $(Y, d_Y, \nu_Y)$ are weight-preserving isometric. The "if" direction is immediate. Note first for the converse that (2.13) implies that for all $\varepsilon > 0$ there exists an $\varepsilon$-isometry from $X$ to $Y$, and therefore, by Lemma 7.3.28 in [11], $d_{\mathrm{GH}}((X, d_X), (Y, d_Y)) < 2\varepsilon$. Thus $d_{\mathrm{GH}}((X, d_X), (Y, d_Y)) = 0$, and it follows from Theorem 7.3.30 of [11] that $(X, d_X)$ and $(Y, d_Y)$ are isometric. Checking the proof of that result, we see that we can construct an isometry $f : X \to Y$ by taking any dense countable set $S \subset X$, any sequence of functions $(f_n)$ such that $f_n$ is an $\varepsilon_n$-isometry with $\varepsilon_n \to 0$ as $n \to \infty$, and letting $f$ be $\lim_k f_{n_k}$ along any subsequence such that the limit exists for all $x \in S$ (such a subsequence exists by the compactness of $Y$). Therefore, fix some dense subset $S \subset X$ and suppose without loss of generality that we have an isometry $f: X \to Y$ given by $f(x) = \lim_{n \to \infty} f_n(x)$, $x \in S$, where $f_n \in F_{X,Y}^{\varepsilon_n}$, $d_P(f_{n*}\nu_X, \nu_Y) \le \varepsilon_n$, and $\lim_{n \to \infty} \varepsilon_n = 0$. We will be done if we can show that $f_*\nu_X = \nu_Y$. If $\mu_X$ is a discrete measure with atoms belonging to $S$, then

$$\begin{aligned}d_P(f_*\nu_X, \nu_Y) &\le \limsup_n [d_P(f_{n*}\nu_X, \nu_Y) + d_P(f_{n*}\mu_X, f_{n*}\nu_X) \\ &\qquad + d_P(f_*\mu_X, f_{n*}\mu_X) + d_P(f_*\nu_X, f_*\mu_X)] \\ &\le 2 d_P(\mu_X, \nu_X),\end{aligned} \quad (2.14)$$

where we have used Lemma 2.2 and the fact that $\lim_{n \to \infty} d_P(f_*\mu_X, f_{n*}\mu_X) = 0$ because of the pointwise convergence of $f_n$ to $f$ on $S$. Because we can choose $\mu_X$ so that $d_P(\mu_X, \nu_X)$ is arbitrarily small, we see that $f_*\nu_X = \nu_Y$, as required.

Now consider three spaces $(X, d_X, \nu_X)$, $(Y, d_Y, \nu_Y)$ and $(Z, d_Z, \nu_Z)$ in $\mathbf{T}^{\mathrm{wt}}$, and constants $\varepsilon, \delta > 0$, such that $\Delta_{\mathrm{GH}^{\mathrm{wt}}}((X, d_X, \nu_X), (Y, d_Y, \nu_Y)) < \varepsilon$ and $\Delta_{\mathrm{GH}^{\mathrm{wt}}}((Y, d_Y, \nu_Y), (Z, d_Z, \nu_Z)) < \delta$. Then there exist $f \in F_{X,Y}^{\varepsilon}$ and $g \in F_{Y,Z}^{\delta}$ such that $d_P(f_*\nu_X, \nu_Y) < \varepsilon$ and $d_P(g_*\nu_Y, \nu_Z) < \delta$. Note that $g \circ f \in F_{X,Z}^{\varepsilon+\delta}$. Moreover, by Lemma 2.2

$$(2.15) \quad d_P((g \circ f)_*\nu_X, \nu_Z) \le d_P(g_*\nu_Y, \nu_Z) + d_P(g_*f_*\nu_X, g_*\nu_Y) < \delta + \varepsilon + \delta.$$

This, and a similar argument with the roles of $X$ and $Z$ interchanged, shows that

$$(2.16) \qquad \Delta_{\mathrm{GH}^{\mathrm{wt}}}(X, Z) \le 2[\Delta_{\mathrm{GH}^{\mathrm{wt}}}(X, Y) + \Delta_{\mathrm{GH}^{\mathrm{wt}}}(Y, Z)].$$

The second inequality in the statement of the lemma is clear. In order to see the first inequality, it suffices to show that for any $Z_1, \ldots, Z_n$ we have

$$(2.17) \qquad \Delta_{\mathrm{GH}^{\mathrm{wt}}}(Z_1, Z_n)^{1/4} \le 2 \sum_{i=1}^{n-1} \Delta_{\mathrm{GH}^{\mathrm{wt}}}(Z_i, Z_{i+1})^{1/4}.$$

We will establish (2.17) by induction. The inequality certainly holds when $n = 2$. Suppose it holds for $2, \ldots, n-1$. Write $S$ for the value of the sum on the right-hand side of (2.17). Put

$$(2.18) \quad k := \max\left\{1 \le m \le n-1 : \sum_{i=1}^{m-1} \Delta_{\mathrm{GH}^{\mathrm{wt}}}(Z_i, Z_{i+1})^{1/4} \le S/2\right\}.$$

By the inductive hypothesis and the definition of $k$,

$$(2.19) \quad \Delta_{\mathrm{GH}^{\mathrm{wt}}}(Z_1, Z_k)^{1/4} \leq 2 \sum_{i=1}^{k-1} \Delta_{\mathrm{GH}^{\mathrm{wt}}}(Z_i, Z_{i+1})^{1/4} \leq 2(S/2) = S.$$

Of course,

$$(2.20) \quad \Delta_{\mathrm{GH}^{\mathrm{wt}}}(Z_k, Z_{k+1})^{1/4} \leq S.$$

By the definition of $k$,

$$(2.21) \quad \sum_{i=1}^{k} \Delta_{\mathrm{GH}^{\mathrm{wt}}}(Z_i, Z_{i+1})^{1/4} > S/2,$$

so that once more by the inductive hypothesis,

$$(2.22) \quad \begin{aligned} \Delta_{\mathrm{GH}^{\mathrm{wt}}}(Z_{k+1}, Z_n)^{1/4} &\leq 2 \sum_{i=k+1}^{n-1} \Delta_{\mathrm{GH}^{\mathrm{wt}}}(Z_i, Z_{i+1})^{1/4} \\ &= 2S - 2 \sum_{i=1}^{k} \Delta_{\mathrm{GH}^{\mathrm{wt}}}(Z_i, Z_{i+1})^{1/4} \\ &\leq S. \end{aligned}$$

From (2.19), (2.20), (2.22) and two applications of (2.16) we have

$$(2.23) \quad \begin{aligned} \Delta_{\mathrm{GH}^{\mathrm{wt}}}(Z_1, Z_n)^{1/4} &\leq \{4[\Delta_{\mathrm{GH}^{\mathrm{wt}}}(Z_1, Z_k) + \Delta_{\mathrm{GH}^{\mathrm{wt}}}(Z_k, Z_{k+1}) \\ &\qquad\qquad + \Delta_{\mathrm{GH}^{\mathrm{wt}}}(Z_{k+1}, Z_n)]\}^{1/4} \\ &\leq (4 \times 3 \times S^4)^{1/4} \\ &\leq 2S, \end{aligned}$$

as required.

It is obvious by construction that $d_{\mathrm{GH}^{\mathrm{wt}}}$ satisfies the triangle inequality. The other properties of a metric follow from the corresponding properties we have already established for $\Delta_{\mathrm{GH}^{\mathrm{wt}}}$ and the bounds in the statement of the lemma which we have already established. □

The procedure we used to construct the weighted Gromov–Hausdorff metric $d_{\mathrm{GH}^{\mathrm{wt}}}$ from the semimetric $\Delta_{\mathrm{GH}^{\mathrm{wt}}}$ was adapted from a proof in [24] of the celebrated result of Alexandroff and Urysohn on the metrizability of uniform spaces. That proof was, in turn, adapted from earlier work of Frink and Bourbaki. The choice of the power $\frac{1}{4}$ is not particularly special; any sufficiently small power would have worked.

Theorem 2.5 below says that the metric space $(\mathbf{T}^{\mathrm{wt}}, d_{\mathrm{GH}^{\mathrm{wt}}})$ is complete and separable and hence is a reasonable space on which to do probability theory. In order to prove this result, we need a compactness criterion that will be useful in its own right.

PROPOSITION 2.4. *A subset $\mathbf{D}$ of $(\mathbf{T}^{\mathrm{wt}}, d_{\mathrm{GH}^{\mathrm{wt}}})$ is relatively compact if and only if the subset $\mathbf{E} := \{(T,d) : (T,d,\nu) \in \mathbf{D}\}$ in $(\mathbf{T}, d_{\mathrm{GH}})$ is relatively compact.*

PROOF. The "only if" direction is clear. Assume for the converse that $\mathbf{E}$ is relatively compact. Suppose that $((T_n, d_{T_n}, \nu_{T_n}))_{n \in \mathbb{N}}$ is a sequence in $\mathbf{D}$. By assumption, $((T_n, d_{T_n}))_{n \in \mathbb{N}}$ has a subsequence converging to some point $(T, d_T)$ of $(\mathbf{T}, d_{\mathrm{GH}})$. For ease of notation, we will renumber and also denote this subsequence by $((T_n, d_{T_n}))_{n \in \mathbb{N}}$. For brevity, we will also omit specific mention of the metric on a real tree when it is clear from the context.

By Proposition 7.4.12 in [11], for each $\varepsilon > 0$ there is a finite $\varepsilon$-net $T^\varepsilon$ in $T$ and for each $n \in \mathbb{N}$ a finite $\varepsilon$-net $T_n^\varepsilon := \{x_n^{\varepsilon,1}, \ldots, x_n^{\varepsilon, \#T_n^\varepsilon}\}$ in $T_n$ such that $d_{\mathrm{GH}}(T_n^\varepsilon, T^\varepsilon) \to 0$ as $n \to \infty$. Moreover, we take $\#T_n^\varepsilon = \#T^\varepsilon = N^\varepsilon$, say, for $n$ sufficiently large, and so, by passing to a further subsequence if necessary, we may assume that $\#T_n^\varepsilon = \#T^\varepsilon = N^\varepsilon$ for all $n \in \mathbb{N}$. We may then assume that $T_n^\varepsilon$ and $T^\varepsilon$ have been indexed so that $\lim_{n \to \infty} d_{T_n}(x_n^{\varepsilon,i}, x_n^{\varepsilon,j}) = d_T(x^{\varepsilon,i}, x^{\varepsilon,j})$ for $1 \le i, j \le N^\varepsilon$.

We may begin with the balls of radius $\varepsilon$ around each point of $T_n^\varepsilon$ and decompose $T_n$ into $N^\varepsilon$ possibly empty, disjoint, measurable sets $\{T_n^{\varepsilon,1}, \ldots, T_n^{\varepsilon, N^\varepsilon}\}$ of radius no greater than $\varepsilon$. Define a measurable map $f_n^\varepsilon : T_n \to T_n^\varepsilon$ by $f_n^\varepsilon(x) = x_n^{\varepsilon,i}$ if $x \in T_n^{\varepsilon,i}$ and let $g_n^\varepsilon$ be the inclusion map from $T_n^\varepsilon$ to $T_n$. By construction, $f_n^\varepsilon$ and $g_n^\varepsilon$ are measurable $\varepsilon$-isometries. Moreover, $d_{\mathrm{P}}((g_n^\varepsilon)_*(f_n^\varepsilon)_*\nu_n, \nu_n) < \varepsilon$ and, of course, $d_{\mathrm{P}}((f_n^\varepsilon)_*\nu_n, (f_n^\varepsilon)_*\nu_n) = 0$. Thus,

$$\Delta_{\mathrm{GH}^{\mathrm{wt}}}((T_n^\varepsilon, (f_n^\varepsilon)_*\nu_n), (T_n, \nu_n)) \le \varepsilon.$$

By similar reasoning, if we define $h_n^\varepsilon : T_n^\varepsilon \to T^\varepsilon$ by $x_n^{\varepsilon,i} \mapsto x^{\varepsilon,i}$, then

$$\lim_{n \to \infty} \Delta_{\mathrm{GH}^{\mathrm{wt}}}((T_n^\varepsilon, (f_n^\varepsilon)_*\nu_n), (T^\varepsilon, (h_n^\varepsilon)_*\nu_n)) = 0.$$

Since $T^\varepsilon$ is finite, by passing to a subsequence (and relabeling as before) we have

$$\lim_{n \to \infty} d_{\mathrm{P}}((h_n^\varepsilon)_*\nu_n, \nu^\varepsilon) = 0$$

for some probability measure $\nu^\varepsilon$ on $T^\varepsilon$, and hence

$$\lim_{n \to \infty} \Delta_{\mathrm{GH}^{\mathrm{wt}}}((T^\varepsilon, (h_n^\varepsilon)_*\nu_n), (T^\varepsilon, \nu^\varepsilon)) = 0.$$

Therefore, by Lemma 2.3,

$$\limsup_{n \to \infty} d_{\mathrm{GH}^{\mathrm{wt}}}((T_n, \nu_n), (T^\varepsilon, (h_n^\varepsilon)_*\nu_n)) \le \varepsilon^{1/4}.$$

Now, since $(T, d_T)$ is compact, the family of measures $\{\nu^\varepsilon : \varepsilon > 0\}$ is relatively compact, and so there is a probability measure $\nu$ on $T$ such that $\nu^\varepsilon$ converges to $\nu$ in the Prohorov distance along a subsequence $\varepsilon \downarrow 0$ and hence, by arguments similar to the above, along the same subsequence $\Delta_{\mathrm{GH^{wt}}}((T^\varepsilon, \nu^\varepsilon), (T, \nu))$ converges to 0. Again applying Lemma 2.3, we have that $d_{\mathrm{GH^{wt}}}((T^\varepsilon, \nu^\varepsilon), (T, \nu))$ converges to 0 along this subsequence.

Combining the foregoing, we see that by passing to a suitable subsequence and relabeling, $d_{\mathrm{GH^{wt}}}((T_n, \nu_n), (T, \nu))$ converges to 0, as required. $\square$

THEOREM 2.5. *The metric space* $(\mathbf{T}^{\mathrm{wt}}, d_{\mathrm{GH^{wt}}})$ *is complete and separable.*

PROOF. Separability follows readily from separability of $(\mathbf{T}, d_{\mathrm{GH}})$ (see Theorem 1 in [20]), and the separability with respect to the Prohorov distance of the probability measures on a fixed complete, separable metric space (see, e.g., [19]), and Lemma 2.3.

It remains to establish completeness. By a standard argument, it suffices to show that any Cauchy sequence in $\mathbf{T}^{\mathrm{wt}}$ has a convergent subsequence. Let $(T_n, d_{T_n}, \nu_n)_{n \in \mathbb{N}}$ be a Cauchy sequence in $\mathbf{T}^{\mathrm{wt}}$. Then $(T_n, d_{T_n})_{n \in \mathbb{N}}$ is a Cauchy sequence in $\mathbf{T}$ by Lemma 2.3. By Theorem 1 in [20] there is a $T \in \mathbf{T}$ such that $d_{\mathrm{GH}}(T_n, T) \to 0$, as $n \to \infty$. In particular, the sequence $(T_n, d_{T_n})_{n \in \mathbb{N}}$ is relatively compact in $\mathbf{T}$, and therefore, by Proposition 2.4, $(T_n, d_{T_n}, \nu_n)_{n \in \mathbb{N}}$ is relatively compact in $\mathbf{T}^{\mathrm{wt}}$. Thus $(T_n, d_{T_n}, \nu_n)_{n \in \mathbb{N}}$ has a convergent subsequence, as required. $\square$

We conclude this section by giving a necessary and sufficient condition for a subset of $(\mathbf{T}, d_{\mathrm{GH}})$ to be relatively compact, and hence, by Proposition 2.4, a necessary and sufficient condition for a subset of $(\mathbf{T}^{\mathrm{wt}}, d_{\mathrm{GH^{wt}}})$ to be relatively compact.

Fix $(T, d) \in \mathbf{T}$ and, as usual, denote the Borel-$\sigma$-algebra on $T$ by $\mathcal{B}(T)$. Let

$$(2.24) \qquad T^o = \bigcup_{a, b \in T} ]a, b[$$

be the *skeleton* of $T$. Observe that if $T' \subset T$ is a dense countable set, then (2.24) holds with $T$ replaced by $T'$. In particular, $T^o \in \mathcal{B}(T)$ and $\mathcal{B}(T)|_{T^o} = \sigma(\{]a, b[; a, b \in T'\})$, where $\mathcal{B}(T)|_{T^o} := \{A \cap T^o; A \in \mathcal{B}(T)\}$. Hence there exists a unique $\sigma$-finite measure $\mu^T$ on $T$, called *length measure*, such that $\mu^T(T \setminus T^o) = 0$ and

$$(2.25) \qquad \mu^T(]a, b[) = d(a, b) \qquad \forall a, b \in T.$$

Such a measure may be constructed as the trace onto $T^o$ of a one-dimensional Hausdorff measure on $T$, and a standard monotone class argument shows that this is the unique measure with property (2.25).

For $\varepsilon > 0$, $T \in \mathbf{T}$ and $\rho \in T$ write

$$(2.26) \qquad R_\varepsilon(T, \rho) := \{x \in T : \exists y \in T, [\rho, y] \ni x, d_T(x, y) \geq \varepsilon\} \cup \{\rho\}$$

4SUBTREE PRUNE AND REGRAFT 13

for the $\varepsilon$-*trimming relative to the root* $\rho$ of the compact $\mathbb{R}$-tree $T$. Then set

$$(2.27) \qquad R_\varepsilon(T) := \begin{cases} \bigcap_{\rho \in T} R_\varepsilon(T, \rho), & \mathrm{diam}(T) > \varepsilon, \\ \mathrm{singleton}, & \mathrm{diam}(T) \leq \varepsilon, \end{cases}$$

where by *singleton* we mean the trivial $\mathbb{R}$-tree consisting of one point. The tree $R_\varepsilon(T)$ is called the $\varepsilon$-*trimming* of the compact $\mathbb{R}$-tree $T$.

LEMMA 2.6. *A subset* $\mathbf{E}$ *of* $(\mathbf{T}, d_{\mathrm{GH}})$ *is relatively compact if and only if for all* $\varepsilon > 0$,

$$(2.28) \qquad \sup\{\mu^T(R_\varepsilon(T)) : T \in \mathbf{E}\} < \infty.$$

PROOF. The "only if" direction follows from the fact that $T \mapsto \mu^T(R_\varepsilon(T))$ is continuous, which is essentially Lemma 7.3 of [20].

Conversely, suppose that (2.28) holds. Given $T \in \mathbf{E}$, an $\varepsilon$-net for $R_\varepsilon(T)$ is a $2\varepsilon$-net for $T$. By Lemma 2.7 below, $R_\varepsilon(T)$ has an $\varepsilon$-net of cardinality at most $[(\frac{\varepsilon}{2})^{-1}\mu^T(R_\varepsilon(T))][(\frac{\varepsilon}{2})^{-1}\mu^T(R_\varepsilon(T)) + 1]$. By assumption, the last quantity is uniformly bounded in $T \in \mathbf{E}$. Thus $\mathbf{E}$ is uniformly totally bounded and hence is relatively compact by Theorem 7.4.15 of [11]. $\square$

LEMMA 2.7. *Let* $T \in \mathbf{T}$ *be such that* $\mu^T(T) < \infty$. *For each* $\varepsilon > 0$ *there is an* $\varepsilon$-*net for* $T$ *of cardinality at most* $[(\frac{\varepsilon}{2})^{-1}\mu^T(T)][(\frac{\varepsilon}{2})^{-1}\mu^T(T) + 1]$.

PROOF. Note that an $\frac{\varepsilon}{2}$-net for $R_{\varepsilon/2}(T)$ will be an $\varepsilon$-net for $T$. The set $T \setminus R_{\varepsilon/2}(T)$ is a collection of disjoint subtrees, one for each leaf of $R_{\varepsilon/2}(T)$, and each such subtree is of diameter at least $\frac{\varepsilon}{2}$. Thus the number of leaves of $R_{\varepsilon/2}(T)$ is at most $(\frac{\varepsilon}{2})^{-1}\mu^T(T)$. Enumerate the leaves of $R_{\varepsilon/2}(T)$ as $x_0, x_1, \ldots, x_n$. Each arc $[x_0, x_i]$, $1 \leq i \leq n$, of $R_{\varepsilon/2}(T)$ has an $\frac{\varepsilon}{2}$-net of cardinality at most $(\frac{\varepsilon}{2})^{-1}d_T(x_0, x_i) + 1 \leq (\frac{\varepsilon}{2})^{-1}\mu^T(T) + 1$. Therefore, by taking the union of these nets, $R_{\varepsilon/2}(T)$ has an $\frac{\varepsilon}{2}$-net of cardinality at most $[(\frac{\varepsilon}{2})^{-1}\mu^T(T)][(\frac{\varepsilon}{2})^{-1}\mu^T(T) + 1]$. $\square$

REMARK 2.8. The bound in Lemma 2.7 is far from optimal. It can be shown that $T$ has an $\varepsilon$-net with a cardinality that is of order $\mu^T(T)/\varepsilon$. This is clear for finite trees (i.e., trees with a finite number of branch points), where we can traverse the tree with a unit speed path and hence think of the tree as an image of the interval $[0, 2\mu^T(T)]$ by a Lipschitz map with Lipschitz constant 1, so that a covering of the interval $[0, 2\mu^T(T)]$ by $\varepsilon$-balls gives a covering of $T$ by $\varepsilon$-balls. This argument can be extended to arbitrary finite length $\mathbb{R}$-trees, but the details are tedious and so we have contented ourselves with the above simpler bound.

**3. Trees and continuous paths.** For the sake of completeness and to establish some notation we recall some facts about the connection between continuous excursion paths and trees (see [4, 18, 26] for more on this connection).

Write $C(\mathbb{R}_+)$ for the space of continuous functions from $\mathbb{R}_+$ into $\mathbb{R}$. For $e \in C(\mathbb{R}_+)$, put $\zeta(e) := \inf\{t > 0 : e(t) = 0\}$ and write

$$
(3.1) \qquad U := \left\{ e \in C(\mathbb{R}_+) : \begin{array}{c} e(0) = 0, \zeta(e) < \infty, \\ e(t) > 0 \text{ for } 0 < t < \zeta(e), \\ \text{and } e(t) = 0 \text{ for } t \geq \zeta(e) \end{array} \right\}
$$

for the space of positive excursion paths. Set $U^\ell := \{e \in U : \zeta(e) = \ell\}$.

We associate each $e \in U^1$ with a compact $\mathbb{R}$-tree as follows. Define an equivalence relation $\sim_e$ on $[0, 1]$ by letting

$$
(3.2) \qquad u_1 \sim_e u_2 \quad \text{iff } e(u_1) = \inf_{u \in [u_1 \wedge u_2, u_1 \vee u_2]} e(u) = e(u_2).
$$

Consider the following pseudometric on $[0, 1]$:

$$
(3.3) \qquad d_{T_e}(u_1, u_2) := e(u_1) - 2 \inf_{u \in [u_1 \wedge u_2, u_1 \vee u_2]} e(u) + e(u_2),
$$

which becomes a true metric on the quotient space $T_e := \mathbb{R}_+|_{\sim_e} = [0, 1]|_{\sim_e}$.

LEMMA 3.1. *For each $e \in U^1$ the metric space $(T_e, d_{T_e})$ is a compact $\mathbb{R}$-tree.*

PROOF. It is straightforward to check that the quotient map from $[0, 1]$ onto $T_e$ is continuous with respect to $d_{T_e}$. Thus $(T_e, d_{T_e})$ is path-connected and compact as the continuous image of a metric space with these properties. In particular, $(T_e, d_{T_e})$ is complete.

To complete the proof, it therefore suffices to verify the *four point condition* (2.1). However, for $u_1, u_2, u_3, u_4 \in T_e$ we have

$$
(3.4) \qquad \begin{aligned} &\max\{d_{T_e}(u_1, u_3) + d_{T_e}(u_2, u_4), d_{T_e}(u_1, u_4) + d_{T_e}(u_2, u_3)\} \\ &\geq d_{T_e}(u_1, u_2) + d_{T_e}(u_3, u_4), \end{aligned}
$$

where strict inequality holds if and only if

$$
(3.5) \qquad \begin{aligned} &\min_{i \neq j} \inf_{u \in [u_i \wedge u_j, u_i \vee u_j]} e(u) \\ &\notin \left\{ \inf_{u \in [u_1 \wedge u_2, u_1 \vee u_2]} e(u), \inf_{u \in [u_3 \wedge u_4, u_3 \vee u_4]} e(u) \right\}. \end{aligned} \qquad \square
$$

REMARK 3.2. Any compact $\mathbb{R}$-tree $T$ is isometric to $T_e$ for some $e \in U^1$. To see this, fix a root $\rho \in T$. Recall $R_\varepsilon(T, \rho)$, the $\varepsilon$-trimming of $T$ with respect to $\rho$ defined in (2.26). Let $\bar{\mu}$ be a probability measure on $T$ that is equivalent to the length measure $\mu^T$. Because $\mu^T$ is $\sigma$-finite, such a probability measure always exists, but one can construct $\bar{\mu}$ explicitly as follows: set $H := \max_{u \in T} d(\rho, u)$, and put

$$\bar{\mu} := 2^{-1} \frac{\mu^T(R(\rho, 2^{-1}H) \cap \cdot)}{\mu^T(R(\rho, 2^{-1}H))}$$
$$+ \sum_{i \geq 2} 2^{-i} \frac{\mu^T(R(\rho, 2^{-i}H) \setminus R(\rho, 2^{-i+1}H) \cap \cdot)}{\mu^T(R(\rho, 2^{-i}H) \setminus R(\rho, 2^{-i+1}H))}.$$

For all $0 < \varepsilon < H$ there is a continuous path

$$f_\varepsilon : [0, 2\mu^T(R_\varepsilon(T, \rho))] \to R_\varepsilon(T, \rho)$$

such that $h_\varepsilon$ defined by $h_\varepsilon(t) := d(\rho, f_\varepsilon(t))$ belongs to $U^{2\mu^T(R_\varepsilon(T,\rho))}$ [in particular, $f_\varepsilon(0) = f_\varepsilon(2\mu^T(R_\varepsilon(T, \rho))) = \rho$], $h_\varepsilon$ is piecewise linear with slopes $\pm 1$ and $T_{h_\varepsilon}$ is isometric to $R_\varepsilon(T, \rho)$. Moreover, these paths may be chosen consistently so that if $\varepsilon' \leq \varepsilon''$, then

$$f_{\varepsilon''}(t) = f_{\varepsilon'}(\inf\{s > 0 : |\{0 \leq r \leq s : f_{\varepsilon'}(r) \in R_{\varepsilon''}(T, \rho)\}| > t\}),$$

where $|\cdot|$ denotes Lebesgue measure. Now define $e_\varepsilon \in U^{\bar{\mu}(R_\varepsilon(T,\rho))}$ to be the absolutely continuous path satisfying

$$\frac{de_\varepsilon(t)}{dt} = 2 \frac{d\mu^T}{d\bar{\mu}}(f_\varepsilon(t)) \frac{dh_\varepsilon(t)}{dt}.$$

It can be shown that $e_\varepsilon$ converges uniformly to some $e \in U^1$ as $\varepsilon \downarrow 0$ and that $T_e$ is isometric to $T$.

From the connection we have recalled between excursion paths and real trees, it should be clear that the analogue of an SPR move for a real tree arising from an excursion path is the excision and reinsertion of a subexcursion. Figure 2 illustrates such an operation.

Each tree coming from a path in $U^1$ has a natural weight on it: for $e \in U^1$, we equip $(T_e, d_{T_e})$ with the weight $\nu_{T_e}$ given by the push-forward of Lebesgue measure on $[0, 1]$ by the quotient map.

We finish this section with a remark about the natural length measure on a tree coming from a path. Given $e \in U^1$ and $a \geq 0$, let

(3.6) $$\mathcal{G}_a := \left\{ t \in [0,1] : \begin{array}{c} e(t) = a \text{ and, for some } \varepsilon > 0, \\ e(u) > a \text{ for all } u \in ]t, t+\varepsilon[, \\ e(t+\varepsilon) = a, \end{array} \right\}$$

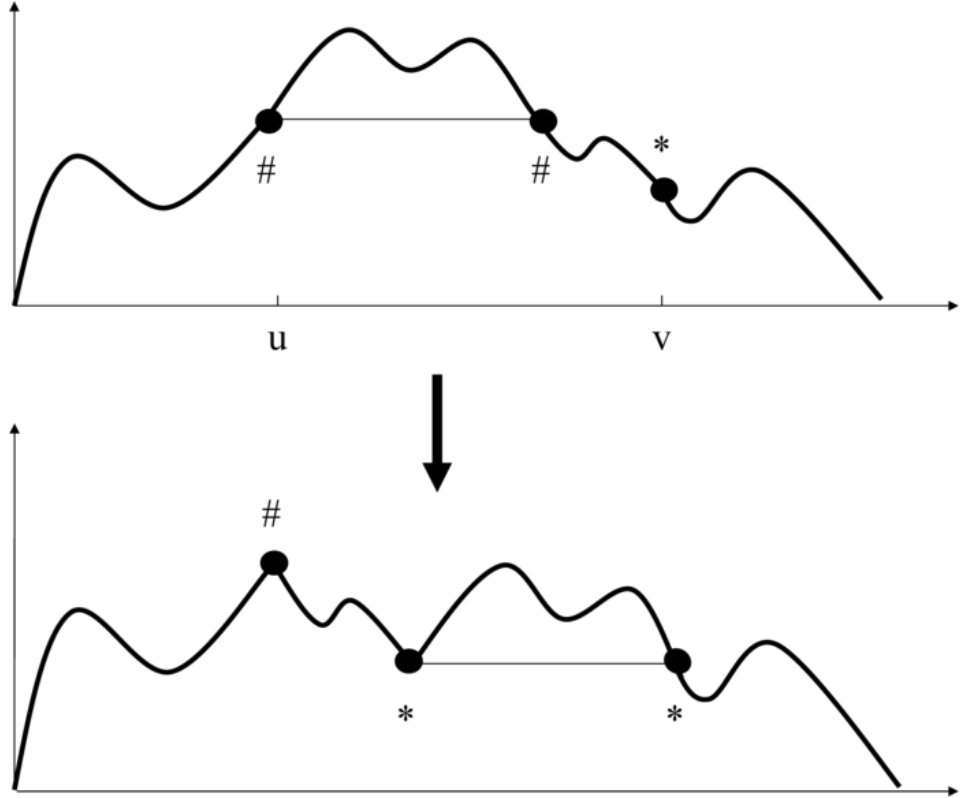

FIG. 2. *A subtree prune and regraft operation on an excursion path: the excursion starting at time u in the top picture is excised and inserted at time v, and the resulting gap between the two points marked # is closed up. The two points marked # (resp. ∗) in the top (resp. bottom) picture correspond to a single point in the associated $\mathbb{R}$-tree.*

denote the countable set of starting points of excursions of the function $e$ above the level $a$. Then $\mu^{T_e}$, the length measure on $T_e$, is just the push-forward of the measure $\int_0^\infty da \sum_{t \in \mathcal{G}_a} \delta_t$ by the quotient map. Alternatively, write

(3.7) $$\Gamma_e := \{(s, a) : s \in \,]0, 1[, a \in [0, e(s)[\,\}$$

for the region between the time axis and the graph of $e$, and for $(s, a) \in \Gamma_e$ denote by $\underline{s}(e, s, a) := \sup\{r < s : e(r) = a\}$ and $\bar{s}(e, s, a) := \inf\{t > s : e(t) = a\}$ the start and finish of the excursion of $e$ above level $a$ that straddles time $s$. Then $\mu^{T_e}$ is the push-forward of the measure $\int_{\Gamma_e} ds \otimes da \, \frac{1}{\bar{s}(e,s,a) - \underline{s}(e,s,a)} \delta_{\underline{s}(e,s,a)}$ by the quotient map. We note that the measure $\mu^{T_e}$ appears in [1].

**4. Uniform random weighted compact $\mathbb{R}$-trees: the continuum random tree.** In this section we will recall the definition of Aldous' continuum ran-dom tree, which can be thought of as a uniformly chosen random weighted compact $\mathbb{R}$-tree.

Consider the Itô excursion measure for excursions of standard Brownian motion away from 0. This $\sigma$-finite measure is defined subject to a normalization of Brownian local time at 0, and we take the usual normalization of local times at each level which makes the local time process an occupation density in the spatial variable for each fixed value of the time variable. The excursion measure is the sum of two measures, one which is concentrated on nonnegative excursions and one which is concentrated on nonpositive excursions. Let $\mathbb{N}$ be the part which is concentrated on nonnegative excursions. Thus, in the notation of Section 3, $\mathbb{N}$ is a $\sigma$-finite measure on $U$, where we equip $U$ with the $\sigma$-field $\mathcal{U}$ generated by the coordinate maps.

Define a map $v : U \to U^1$ by $e \mapsto \frac{e(\zeta(e)\cdot)}{\sqrt{\zeta(e)}}$. Then

$$(4.1) \qquad \mathbb{P}(\Gamma) := \frac{\mathbb{N}\{v^{-1}(\Gamma) \cap \{e \in U : \zeta(e) \geq c\}\}}{\mathbb{N}\{e \in U : \zeta(e) \geq c\}}, \qquad \Gamma \in \mathcal{U}_{|U^1},$$

does not depend on $c > 0$ (see, e.g., Exercise 12.2.13.2 in [27]). The probability measure $\mathbb{P}$ is called the law of normalized nonnegative Brownian excursion. We have

$$(4.2) \qquad \mathbb{N}\{e \in U : \zeta(e) \in dc\} = \frac{dc}{2\sqrt{2\pi c^3}}$$

and, defining $\mathcal{S}_c : U^1 \to U^c$ by

$$(4.3) \qquad \mathcal{S}_c e := \sqrt{c} e(\cdot / c),$$

we have

$$(4.4) \qquad \int \mathbb{N}(de) G(e) = \int_0^\infty \frac{dc}{2\sqrt{2\pi c^3}} \int_{U^1} \mathbb{P}(de) G(\mathcal{S}_c e)$$

for a nonnegative measurable function $G : U \to \mathbb{R}$.

Recall from Section 3 how each $e \in U^1$ is associated with a weighted compact $\mathbb{R}$-tree $(T_e, d_{T_e}, \nu_{T_e})$. Let $\mathbf{P}$ be the probability measure on $(\mathbf{T}^{\mathrm{wt}}, d_{\mathrm{GH}^{\mathrm{wt}}})$ that is the push-forward of the normalized excursion measure by the map $e \mapsto (T_{2e}, d_{T_{2e}}, \nu_{T_{2e}})$, where $2e \in U^1$ is just the excursion path $t \mapsto 2e(t)$.

The probability measure $\mathbf{P}$ is the distribution of an object consisting of Aldous' *continuum random tree* along with a natural measure on this tree (see, e.g., [2, 4]). The continuum random tree arises as the limit of a uniform random tree on $n$ vertices when $n \to \infty$ and edge lengths are rescaled by a factor of $1/\sqrt{n}$. The appearance of $2e$ rather than $e$ in the definition of $\mathbf{P}$ is a consequence of this choice of scaling. The associated probability measure on each realization of the continuum random tree is the measure that arises in this limiting construction by taking the uniform probability measure on realizations of the approximating finite trees. The probability measure $\mathbf{P}$ can therefore be viewed informally as the "uniform distribution" on $(\mathbf{T}^{\mathrm{wt}}, d_{\mathrm{GH}^{\mathrm{wt}}})$.

**5. Campbell measure facts.** For the purposes of constructing the Markov process that is of interest to us, we need to understand picking a random weighted tree $(T, d_T, \nu_T)$ according to the continuum random tree distribution $\mathbf{P}$, picking a point $u$ according to the length measure $\mu^T$ and another point $v$ according to the weight $\nu_T$, and then decomposing $T$ into two subtrees rooted at $u$—one that contains $v$ and one that does not (we are being a little imprecise here, because $\mu^T$ will be an infinite measure, $\mathbf{P}$ almost surely).

In order to understand this decomposition, we must understand the corresponding decomposition of excursion paths under normalized excursion measure. Because subtrees correspond to subexcursions and because of our observation in Section 3 that for an excursion $e$ the length measure $\mu^{T_e}$ on the corresponding tree is the push-forward of the measure $\int_{\Gamma_e} ds \otimes da \, \frac{1}{\bar{s}(e,s,a) - \underline{s}(e,s,a)} \delta_{\underline{s}(e,s,a)}$ by the quotient map, we need to understand the decomposition of the excursion $e$ into the excursion above $a$ that straddles $s$ and the "remaining" excursion when $e$ is chosen according to the standard Brownian excursion distribution $\mathbb{P}$ and $(s,a)$ is chosen according to the $\sigma$-finite measure $ds \otimes da \, \frac{1}{\bar{s}(e,s,a) - \underline{s}(e,s,a)}$ on $\Gamma_e$—see Figure 3.

Given an excursion $e \in U$ and a level $a \geq 0$ write:

(a) $\zeta(e) := \inf\{t > 0 : e(t) = 0\}$ for the "length" of $e$,
(b) $\ell_t^a(e)$ for the local time of $e$ at level $a$ up to time $t$,
(c) $e^{\downarrow a}$ for $e$ time-changed by the inverse of $t \mapsto \int_0^t ds \, \mathbb{1}\{e(s) \leq a\}$ (i.e., $e^{\downarrow a}$ is $e$ with the subexcursions above level $a$ excised and the gaps closed up),
(d) $\ell_t^a(e^{\downarrow a})$ for the local time of $e^{\downarrow a}$ at the level $a$ up to time $t$,
(e) $U^{\uparrow a}(e)$ for the set of subexcursion intervals of $e$ above $a$ (i.e., an element of $U^{\uparrow a}(e)$ is an interval $I = [g_I, d_I]$ such that $e(g_I) = e(d_I) = a$ and $e(t) > a$ for $g_I < t < d_I$),
(f) $\mathcal{N}^{\uparrow a}(e)$ for the counting measure that puts a unit mass at each point $(s', e')$, where, for some $I \in U^{\uparrow a}(e)$, $s' := \ell_{g_I}^a(e)$ is the amount of local time of $e$ at level $a$ accumulated up to the beginning of the subexcursion $I$ and $e' \in U$ given by

$$(5.1) \qquad e'(t) = \begin{cases} e(g_I + t) - a, & 0 \leq t \leq d_I - g_I, \\ 0, & t > d_I - g_I, \end{cases}$$

is the corresponding piece of the path $e$ shifted to become an excursion above the level 0 starting at time 0,

(g) $\hat{e}^{s,a} \in U$ and $\check{e}^{s,a} \in U$, for the subexcursion "above" $(s,a) \in \Gamma_e$, that is,

$$(5.2) \quad \hat{e}^{s,a}(t) := \begin{cases} e(\underline{s}(e,s,a) + t) - a, & 0 \leq t \leq \bar{s}(e,s,a) - \underline{s}(e,s,a), \\ 0, & t > \bar{s}(e,s,a) - \underline{s}(e,s,a), \end{cases}$$

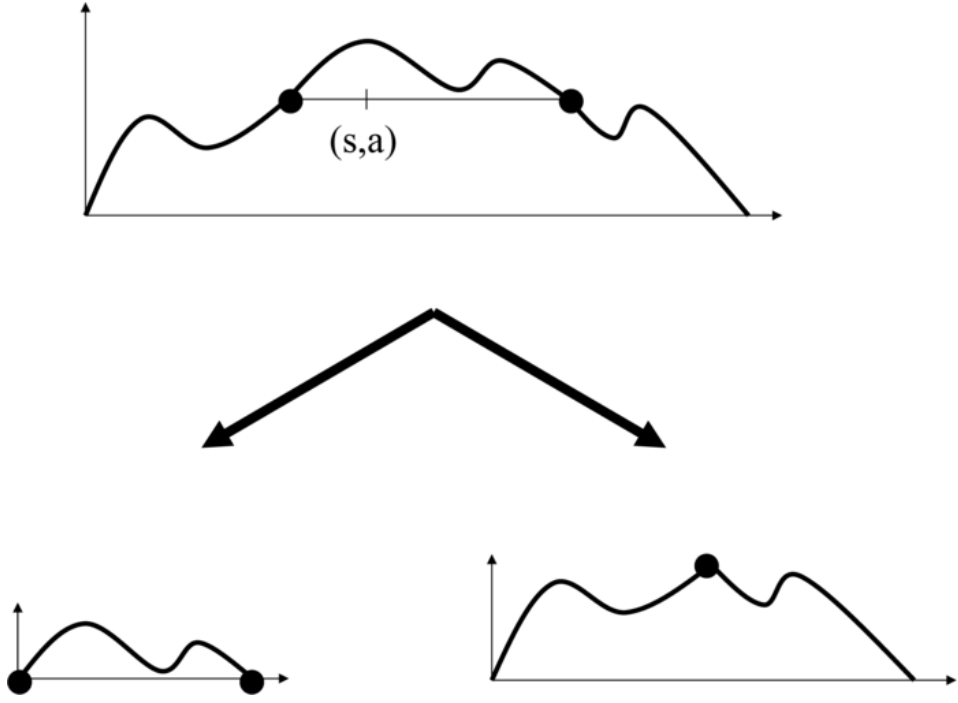

FIG. 3. *The decomposition of the excursion e (top picture) into the excursion $\hat{e}^{s,a}$ above level a that straddles time s (bottom left picture) and the "remaining" excursion $\check{e}^{s,a}$ (bottom right picture).*

respectively "below" $(s,a) \in \Gamma_e$, that is,

$$(5.3) \quad \check{e}^{s,a}(t) := \begin{cases} e(t), & 0 \le t \le \underline{s}(e,s,a), \\ e(t + \bar{s}(e,s,a) - \underline{s}(e,s,a)), & t > \underline{s}(e,s,a), \end{cases}$$

(h) $\sigma_s^a(e) := \inf\{t \ge 0 : \ell_t^a(e) \ge s\}$ and $\tau_s^a(e) := \inf\{t \ge 0 : \ell_t^a(e) > s\}$,

(i) $\tilde{e}^{s,a} \in U$ for $e$ with the interval $]\sigma_s^a(e), \tau_s^a(e)[$ containing an excursion above level $a$ excised, that is,

$$(5.4) \quad \tilde{e}^{s,a}(t) := \begin{cases} e(t), & 0 \le t \le \sigma_s^a(e), \\ e(t + \tau_s^a(e) - \sigma_s^a(e)), & t > \sigma_s^a(e). \end{cases}$$

The following path decomposition result under the $\sigma$-finite measure $\mathbb{N}$ is preparatory to a decomposition under the probability measure $\mathbb{P}$, Corollary 5.2, that has a simpler intuitive interpretation.

PROPOSITION 5.1. *For nonnegative measurable functions $F$ on $\mathbb{R}_+$ and $G, H$ on $U$,*

$$\int \mathbb{N}(de) \int_{\Gamma_e} \frac{ds \otimes da}{\bar{s}(e,s,a) - \underline{s}(e,s,a)} F(\underline{s}(e,s,a)) G(\hat{e}^{s,a}) H(\check{e}^{s,a})$$

$$= \int \mathbb{N}(de) \int_0^\infty da \int \mathcal{N}^{\uparrow a}(e)(d(s',e'))F(\sigma_{s'}^a(e))G(e')H(\tilde{e}^{s',a})$$

$$= \mathbb{N}[G]\mathbb{N}\left[H \int_0^\zeta ds F(s)\right].$$

PROOF. The first equality is just a change in the order of integration and has already been remarked upon in Section 3.

Standard excursion theory (see, e.g., [8, 27, 28]) says that under $\mathbb{N}$, the random measure $e \mapsto \mathcal{N}^{\uparrow a}(e)$ conditional on $e \mapsto e^{\downarrow a}$ is a Poisson random measure with intensity measure $\lambda^{\downarrow a}(e) \otimes \mathbb{N}$, where $\lambda^{\downarrow a}(e)$ is Lebesgue measure restricted to the interval $[0, \ell_\infty^a(e)] = [0, 2\ell_\infty^a(e^{\downarrow a})]$.

Note that $\tilde{e}^{s',a}$ is constructed from $e^{\downarrow a}$ and $\mathcal{N}^{\uparrow a}(e) - \delta_{(s',e')}$ in the same way that $e$ is constructed from $e^{\downarrow a}$ and $\mathcal{N}^{\uparrow a}(e)$. Also, $\sigma_{s'}^a(\tilde{e}^{s',a}) = \sigma_{s'}^a(e)$. Therefore, by the Campbell–Palm formula for Poisson random measures (see, e.g., Section 12.1 of [13]),

$$\int \mathbb{N}(de) \int_0^\infty da \int \mathcal{N}^{\uparrow a}(e)(d(s',e'))F(\sigma_{s'}^a(e))G(e')H(\tilde{e}^{s',a})$$

$$= \int \mathbb{N}(de) \int_0^\infty da \, \mathbb{N}\left[\int \mathcal{N}^{\uparrow a}(e)(d(s',e'))F(\sigma_{s'}^a(e))G(e')H(\tilde{e}^{s',a})\Big|e^{\downarrow a}\right]$$

$$= \int \mathbb{N}(de) \int_0^\infty da \, \mathbb{N}[G]\mathbb{N}\left[\left\{\int_0^{\ell_\infty^a(e)} ds' F(\sigma_{s'}^a(e))\right\}H\Big|e^{\downarrow a}\right]$$

$$= \mathbb{N}[G] \int_0^\infty da \int \mathbb{N}(de)\left(\left\{\int d\ell_s^a(e)F(s)\right\}H(e)\right)$$

$$= \mathbb{N}[G] \int \mathbb{N}(de)\left(\left\{\int_0^\infty da \int d\ell_s^a(e)F(s)\right\}H(e)\right)$$

$$= \mathbb{N}[G]\mathbb{N}\left[H \int_0^\zeta ds F(s)\right]. \qquad \square$$

The next result says that if we pick an excursion $e$ according to the standard excursion distribution $\mathbb{P}$ and then pick a point $(s,a) \in \Gamma_e$ according to the $\sigma$-finite length measure corresponding to the length measure $\mu^{T_e}$ on the associated tree $T_e$ (see the end of Section 3), then the following objects are independent:

(a) the length of the excursion above level $a$ that straddles time $s$,

(b) the excursion obtained by taking the excursion above level $a$ that straddles time $s$, turning it (by a shift of axes) into an excursion $\hat{e}^{s,a}$ above level zero starting at time zero, and then Brownian rescaling $\hat{e}^{s,a}$ to produce an excursion of unit length,

(c) the excursion obtained by taking the excursion $\check{e}^{s,a}$ that comes from excising $\hat{e}^{s,a}$ and closing up the gap, and then Brownian rescaling $\check{e}^{s,a}$ to produce an excursion of unit length,

(d) the starting time $\underline{s}(e, s, a)$ of the excursion above level $a$ that straddles time $s$ rescaled by the length of $\check{e}^{s,a}$ to give a time in the interval $[0, 1]$.

Moreover, the length in (a) is "distributed" according to the $\sigma$-finite measure

$$(5.5) \qquad \frac{1}{2\sqrt{2\pi}} \frac{d\rho}{\sqrt{(1-\rho)\rho^3}}, \qquad 0 \leq \rho \leq 1,$$

the unit length excursions in (b) and (c) are both distributed as standard Brownian excursions (i.e., according to $\mathbb{P}$) and the time in (d) is uniformly distributed on the interval $[0, 1]$.

Recall from (4.3) that $\mathcal{S}_c : U^1 \to U^c$ is the Brownian rescaling map defined by

$$\mathcal{S}_c e := \sqrt{c} e(\cdot/c).$$

COROLLARY 5.2. *For nonnegative measurable functions $F$ on $\mathbb{R}_+$ and $K$ on $U \times U$,*

$$\int \mathbb{P}(de) \int_{\Gamma_e} \frac{ds \otimes da}{\bar{s}(e, s, a) - \underline{s}(e, s, a)} F\left(\frac{\underline{s}(e, s, a)}{\zeta(\check{e}^{s,a})}\right) K(\hat{e}^{s,a}, \check{e}^{s,a})$$

$$= \left\{\int_0^1 du F(u)\right\} \int \mathbb{P}(de) \int_{\Gamma_e} \frac{ds \otimes da}{\bar{s}(e, s, a) - \underline{s}(e, s, a)} K(\hat{e}^{s,a}, \check{e}^{s,a})$$

$$= \left\{\int_0^1 du F(u)\right\} \frac{1}{2\sqrt{2\pi}} \int_0^1 \frac{d\rho}{\sqrt{(1-\rho)\rho^3}} \int \mathbb{P}(de') \otimes \mathbb{P}(de'') K(\mathcal{S}_\rho e', \mathcal{S}_{1-\rho} e'').$$

PROOF. For a nonnegative measurable function $L$ on $U \times U$, it follows straightforwardly from Proposition 5.1 that

$$(5.6) \quad \begin{aligned} &\int \mathbb{N}(de) \int_{\Gamma_e} \frac{ds \otimes da}{\bar{s}(e, s, a) - \underline{s}(e, s, a)} F\left(\frac{\underline{s}(e, s, a)}{\zeta(\check{e}^{s,a})}\right) L(\hat{e}^{s,a}, \check{e}^{s,a}) \\ &= \left\{\int_0^1 du\, F(u)\right\} \int \mathbb{N}(de') \otimes \mathbb{N}(de'') L(e', e'') \zeta(e''). \end{aligned}$$

The left-hand side of (5.6) is, by (4.4),

$$(5.7) \quad \begin{aligned} &\int_0^\infty \frac{dc}{2\sqrt{2\pi c^3}} \int \mathbb{P}(de) \int_{\Gamma_{\mathcal{S}_c e}} ds \otimes da \\ &\qquad \times \frac{F(\underline{s}(\mathcal{S}_c e, s, a)/\zeta(\widetilde{\mathcal{S}_c e}^{s,a})) L(\widehat{\mathcal{S}_c e}^{s,a}, \widetilde{\mathcal{S}_c e}^{s,a})}{\bar{s}(\mathcal{S}_c e, s, a) - \underline{s}(\mathcal{S}_c e, s, a)}. \end{aligned}$$

If we change variables to $t = s/c$ and $b = a/\sqrt{c}$, then the integral for $(s, a)$ over $\Gamma_{\mathcal{S}_c e}$ becomes an integral for $(t, b)$ over $\Gamma_e$. Also,

$$\underline{s}(\mathcal{S}_c e, ct, \sqrt{c}b) = \sup\left\{r < ct : \sqrt{c}e\left(\frac{r}{c}\right) < \sqrt{c}b\right\}$$
(5.8)
$$= c \sup\{r < t : e(r) < b\}$$
$$= c\underline{s}(e, t, b),$$

and, by similar reasoning,

(5.9) $$\bar{s}(\mathcal{S}_c e, ct, \sqrt{c}b) = c\bar{s}(e, t, b)$$

and

(5.10) $$\zeta(\widetilde{\mathcal{S}_c e}^{ct, \sqrt{c}b}) = c\zeta(\check{e}^{t,b}).$$

Thus (5.7) is

(5.11)
$$\int_0^\infty \frac{dc}{2\sqrt{2\pi c^3}} \int \mathbb{P}(de) \sqrt{c} \int_{\Gamma_e} dt \otimes db$$
$$\times \frac{F(\underline{s}(e,t,b)/\zeta(\check{e}^{t,b})) L(\widehat{\mathcal{S}_c e}^{ct,\sqrt{c}b}, \widetilde{\mathcal{S}_c e}^{ct,\sqrt{c}b})}{\bar{s}(e,t,b) - \underline{s}(e,t,b)}.$$

Now suppose that $L$ is of the form

(5.12) $$L(e', e'') = K(\mathcal{R}_{\zeta(e')+\zeta(e'')} e', \mathcal{R}_{\zeta(e')+\zeta(e'')} e'') \frac{M(\zeta(e') + \zeta(e''))}{\sqrt{\zeta(e') + \zeta(e'')}},$$

where, for ease of notation, we put for $e \in U$, and $c > 0$,

(5.13) $$\mathcal{R}_c e := \mathcal{S}_{c^{-1}} e = \frac{1}{\sqrt{c}} e(c \cdot).$$

Then (5.11) becomes

(5.14)
$$\int_0^\infty \frac{dc}{2\sqrt{2\pi c^3}} \int \mathbb{P}(de) \int_{\Gamma_e} dt \otimes db$$
$$\times \frac{F(\underline{s}(e,t,b)/\zeta(\check{e}^{t,b})) K(\hat{e}^{t,b}, \check{e}^{t,b}) M(c)}{\bar{s}(e,t,b) - \underline{s}(e,t,b)}.$$

Since (5.14) was shown to be equivalent to the left-hand side of (5.6), it follows from (4.4) that

(5.15)
$$\int \mathbb{P}(de) \int_{\Gamma_e} \frac{dt \otimes db}{\bar{s}(e,t,b) - \underline{s}(e,t,b)} F\left(\frac{\underline{s}(e,t,b)}{\zeta(\check{e}^{t,b})}\right) K(\hat{e}^{t,b}, \check{e}^{t,b})$$
$$= \frac{\int_0^1 du F(u)}{\mathbb{N}[M]} \int \mathbb{N}(de') \otimes \mathbb{N}(de'') L(e', e'') \zeta(e''),$$

and the first equality of the statement follows.

We have from identity (5.15) that, for any $C > 0$,

$$\mathbb{N}\{\zeta(e) > C\} \int \mathbb{P}(de) \int_{\Gamma_e} \frac{ds \otimes da}{\overline{s}(e,s,a) - \underline{s}(e,s,a)} K(\hat{e}^{s,a}, \check{e}^{s,a})$$

$$= \int \mathbb{N}(de') \otimes \mathbb{N}(de'') K(\mathcal{R}_{\zeta(e')+\zeta(e'')}e', \mathcal{R}_{\zeta(e')+\zeta(e'')}e'')$$

$$\times \frac{\mathbb{1}\{\zeta(e') + \zeta(e'') > C\}}{\sqrt{\zeta(e') + \zeta(e'')}} \zeta(e'')$$

$$= \int_0^\infty \frac{dc'}{2\sqrt{2\pi c'^3}} \int_0^\infty \frac{dc''}{2\sqrt{2\pi c''}}$$

$$\times \int \mathbb{P}(de') \otimes \mathbb{P}(de'') K(\mathcal{R}_{c'+c''}\mathcal{S}_{c'}e', \mathcal{R}_{c'+c''}\mathcal{S}_{c''}e'') \frac{\mathbb{1}\{c' + c'' > C\}}{\sqrt{c' + c''}}.$$

Make the change of variables $\rho = \frac{c'}{c'+c''}$ and $\xi = c' + c''$ (with corresponding Jacobian factor $\xi$) to get

$$\int_0^\infty \frac{dc'}{2\sqrt{2\pi c'^3}} \int_0^\infty \frac{dc''}{2\sqrt{2\pi c''}}$$

$$\times \int \mathbb{P}(de') \otimes \mathbb{P}(de'') K(\mathcal{R}_{c'+c''}\mathcal{S}_{c'}e', \mathcal{R}_{c'+c''}\mathcal{S}_{c''}e'') \frac{\mathbb{1}\{c' + c'' > C\}}{\sqrt{c' + c''}}$$

$$= \left(\frac{1}{2\sqrt{2\pi}}\right)^2 \int_0^\infty d\xi \int_0^1 \frac{d\rho\,\xi}{\sqrt{\rho^3(1-\rho)\xi^4}} \frac{\mathbb{1}\{\xi > C\}}{\sqrt{\xi}}$$

$$\times \int \mathbb{P}(de') \otimes \mathbb{P}(de'') K(\mathcal{S}_\rho e', \mathcal{S}_{1-\rho} e'')$$

$$= \left(\frac{1}{2\sqrt{2\pi}}\right)^2 \left\{\int_C^\infty \frac{d\xi}{\sqrt{\xi^3}}\right\} \int_0^1 \frac{d\rho}{\sqrt{\rho^3(1-\rho)}}$$

$$\times \int \mathbb{P}(de') \otimes \mathbb{P}(de'') K(\mathcal{S}_\rho e', \mathcal{S}_{1-\rho} e''),$$

and the corollary follows upon recalling (4.2). $\square$

COROLLARY 5.3. (i) *For* $x > 0$,

$$\int \mathbb{P}(de) \int_{\Gamma_e} \frac{ds \otimes da}{\overline{s}(e,s,a) - \underline{s}(e,s,a)} \mathbb{1}\left\{\max_{0 \le t \le \zeta(\hat{e}^{s,a})} \hat{e}^{s,a} > x\right\}$$

$$= 2 \sum_{n=1}^\infty nx \exp(-2n^2 x^2).$$

(ii) *For $0 < p \leq 1$,*

$$\int \mathbb{P}(de) \int_{\Gamma_e} \frac{ds \otimes da}{\bar{s}(e,s,a) - \underline{s}(e,s,a)} \mathbb{1}\{\zeta(\hat{e}^{s,a}) > p\} = \sqrt{\frac{1-p}{2\pi p}}.$$

PROOF. (i) Recall first of all from Theorem 5.2.10 in [25] that

(5.16) $$\mathbb{P}\left\{e \in U^1 : \max_{0 \leq t \leq 1} e(t) > x\right\} = 2\sum_{n=1}^{\infty}(4n^2 x^2 - 1)\exp(-2n^2 x^2).$$

By Corollary 5.2 applied to $K(e', e'') := \mathbb{1}\{\max_{t \in [0,\zeta(e')]} e'(t) \geq x\}$ and $F \equiv 1$,

$$\int \mathbb{P}(de) \int_{\Gamma_e} \frac{ds \otimes da}{\bar{s}(e,s,a) - \underline{s}(e,s,a)} \mathbb{1}\left\{\max_{0 \leq t \leq \zeta(\hat{e}^{s,a})} \hat{e}^{s,a} > x\right\}$$

$$= \frac{1}{2\sqrt{2\pi}} \int_0^1 \frac{d\rho}{\sqrt{\rho^3(1-\rho)}} \mathbb{P}\left\{\max_{t \in [0,\rho]} \sqrt{\rho} e(t/\rho) > x\right\}$$

$$= \frac{1}{2\sqrt{2\pi}} \int_0^1 \frac{d\rho}{\sqrt{\rho^3(1-\rho)}} \mathbb{P}\left\{\max_{t \in [0,1]} e(t) > \frac{x}{\sqrt{\rho}}\right\}$$

$$= \frac{1}{2\sqrt{2\pi}} \int_0^1 \frac{d\rho}{\sqrt{\rho^3(1-\rho)}} 2 \sum_{n=1}^{\infty}\left(4n^2 \frac{x^2}{\rho} - 1\right)\exp\left(-2n^2 \frac{x^2}{\rho}\right)$$

$$= 2\sum_{n=1}^{\infty} nx \exp(-2n^2 x^2),$$

as claimed.

(ii) Corollary 5.2 applied to $K(e', e'') := \mathbb{1}\{\zeta(e') \geq p\}$ and $F \equiv 1$ immediately yields

$$\int \mathbb{P}(de) \int_{\Gamma_e} \frac{ds \otimes da}{\bar{s}(e,s,a) - \underline{s}(e,s,a)} \mathbb{1}\{\zeta(\hat{e}^{s,a}) > p\}$$

$$= \frac{1}{2\sqrt{2\pi}} \int_p^1 \frac{d\rho}{\sqrt{\rho^3(1-\rho)}} = \sqrt{\frac{1-p}{2\pi p}}. \qquad \square$$

We conclude this section by calculating the expectations of some functionals with respect to **P** [the "uniform distribution" on $(\mathbf{T}^{\mathrm{wt}}, d_{\mathrm{GH^{wt}}})$ as introduced in the end of Section 4].

For $T \in \mathbf{T}^{\mathrm{wt}}$, and $\rho \in T$, recall $R_c(T, \rho)$ from (2.26), and the length measure $\mu^T$ from (2.25). Given $(T, d) \in \mathbf{T}^{\mathrm{wt}}$ and $u, v \in T$, let

(5.17) $$S^{T,u,v} := \{w \in T : u \in \,]v, w[\,\}$$

denote the subtree of $T$ that differs from its closure by the point $u$, which can be thought of as its root, and consists of points that are on the "other side" of $u$ from $v$ (recall $]v, w[$ is the open arc in $T$ between $v$ and $w$).

LEMMA 5.4. (i) *For $x > 0$,*

$$\mathbf{P}[\mu^T \otimes \nu_T \{(u,v) \in T \times T : \text{height}(S^{T,u,v}) > x\}]$$
$$= \mathbf{P}\left[\int_T \nu_T(dv)\mu^T(R_x(T,v))\right]$$
$$= 2\sum_{n=1}^{\infty} nx \exp(-n^2 x^2/2).$$

(ii) *For $1 < \alpha < \infty$,*

$$\mathbf{P}\left[\int_T \nu_T(dv) \int_T \mu^T(du)(\text{height}(S^{T,u,v}))^\alpha\right]$$
$$= 2^{(\alpha+1)/2} \alpha \Gamma\left(\frac{\alpha+1}{2}\right) \zeta(\alpha),$$

*where, as usual, $\zeta(\alpha) := \sum_{n \geq 1} n^{-\alpha}$.*

(iii) *For $0 < p \leq 1$,*

$$\mathbf{P}[\mu^T \otimes \nu_T \{(u,v) \in T \times T : \nu_T(S^{T,u,v}) > p\}] = \sqrt{\frac{2(1-p)}{\pi p}}.$$

(iv) *For $\frac{1}{2} < \beta < \infty$,*

$$\mathbf{P}\left[\int_T \nu_T(dv) \int_T \mu^T(du)(\nu_T(S^{T,u,v}))^\beta\right] = 2^{-1/2} \frac{\Gamma(\beta - 1/2)}{\Gamma(\beta)}.$$

PROOF. (i) The first equality is clear from the definition of $R_x(T,v)$ and Fubini's theorem.

Turning to the equality of the first and last terms, first recall that $\mathbf{P}$ is the push-forward on $(\mathbf{T}^{\text{wt}}, d_{\text{GH}^{\text{wt}}})$ of the normalized excursion measure $\mathbb{P}$ by the map $e \mapsto (T_{2e}, d_{T_{2e}}, \nu_{T_{2e}})$, where $2e \in U^1$ is just the excursion path $t \mapsto 2e(t)$. In particular, $T_{2e}$ is the quotient of the interval $[0,1]$ by the equivalence relation defined by $2e$. By the invariance of the standard Brownian excursion under random rerooting (see Section 2.7 of [3]), the point in $T_{2e}$ that corresponds to the equivalence class of $0 \in [0,1]$ is distributed according to $\nu_{T_{2e}}$ when $e$ is chosen according to $\mathbb{P}$. Moreover, recall from the end of Section 3 that for $e \in U^1$, the length measure $\mu^{T_e}$ is the push-forward of the measure $ds \otimes da \frac{1}{\bar{s}(e,s,a) - \underline{s}(e,s,a)} \delta_{\underline{s}(e,s,a)}$ on the subgraph $\Gamma_e$ by the quotient map defined in (3.2).

It follows that if we pick $T$ according to $\mathbf{P}$ and then pick $(u,v) \in T \times T$ according to $\mu^T \otimes \nu_T$, then the subtree $S^{T,u,v}$ that arises has the same $\sigma$-finite law as the tree associated with the excursion $2\hat{e}^{s,a}$ when $e$ is chosen according to $\mathbb{P}$ and $(s,a)$ is chosen according to the measure $ds \otimes da \frac{1}{\bar{s}(e,s,a) - \underline{s}(e,s,a)} \delta_{\underline{s}(e,s,a)}$ on the subgraph $\Gamma_e$.

Therefore, by part (i) of Corollary 5.3,

$$\mathbf{P}\left[\int_T \nu_T(dv) \int_T \mu^T(du) \mathbb{1}\{\text{height}(S^{T,u,v}) > x\}\right]$$

$$= 2\int \mathbb{P}(de) \int_{\Gamma_e} \frac{ds \otimes da}{\overline{s}(e,s,a) - \underline{s}(e,s,a)} \mathbb{1}\left\{\max_{0 \leq t \leq \zeta(\hat{e}^{s,a})} \hat{e}^{s,a} > \frac{x}{2}\right\}$$

$$= 2\sum_{n=1}^{\infty} nx \exp(-n^2 x^2/2).$$

Part (ii) is a consequence of part (i) and some straightforward calculus.

Part (iii) follows immediately from part (ii) of Corollary 5.3.

Part (iv) is a consequence of part (iii) and some more straightforward calculus. □

**6. A symmetric jump measure on $(\mathbf{T}^{\text{wt}}, d_{\mathbf{GH}^{\text{wt}}})$.** In this section we will construct and study a measure on $\mathbf{T}^{\text{wt}} \times \mathbf{T}^{\text{wt}}$ that is related to the decomposition discussed at the beginning of Section 5.

Define a map $\Theta$ from $\{((T,d),u,v) : T \in \mathbf{T}, u \in T, v \in T\}$ into $\mathbf{T}$ by setting $\Theta((T,d),u,v) := (T,d^{(u,v)})$, letting

$$(6.1) \quad d^{(u,v)}(x,y) := \begin{cases} d(x,y), & \text{if } x,y \in S^{T,u,v}, \\ d(x,y), & \text{if } x,y \in T \setminus S^{T,u,v}, \\ d(x,u) + d(v,y), & \text{if } x \in S^{T,u,v}, y \in T \setminus S^{T,u,v}, \\ d(y,u) + d(v,x), & \text{if } y \in S^{T,u,v}, x \in T \setminus S^{T,u,v}. \end{cases}$$

That is, $\Theta((T,d),u,v)$ is just $T$ as a set, but the metric has been changed so that the subtree $S^{T,u,v}$ with root $u$ is now pruned and regrafted so as to have root $v$.

If $(T,d,\nu) \in \mathbf{T}^{\text{wt}}$ and $(u,v) \in T \times T$, then we can think of $\nu$ as a weight on $(T,d^{(u,v)})$, because the Borel structures induced by $d$ and $d^{(u,v)}$ are the same. With a slight misuse of notation we will therefore write $\Theta((T,d,\nu),u,v)$ for $(T,d^{(u,v)},\nu) \in \mathbf{T}^{\text{wt}}$. Intuitively, the mass contained in $S^{T,u,v}$ is transported along with the subtree.

Define a kernel $\kappa$ on $\mathbf{T}^{\text{wt}}$ by

$$(6.2) \quad \kappa((T,d_T,\nu_T),\mathbf{B}) := \mu^T \otimes \nu_T\{(u,v) \in T \times T : \Theta(T,u,v) \in \mathbf{B}\}$$

for $\mathbf{B} \in \mathcal{B}(\mathbf{T}^{\text{wt}})$. Thus $\kappa((T,d_T,\nu_T),\cdot)$ is the jump kernel described informally in the Introduction.

REMARK 6.1. It is clear that $\kappa((T,d_T,\nu_T),\cdot)$ is a Borel measure on $\mathbf{T}^{\text{wt}}$ for each $(T,d_T,\nu_T) \in \mathbf{T}^{\text{wt}}$. In order to show that $\kappa(\cdot,\mathbf{B})$ is a Borel function on $\mathbf{T}^{\mathrm{wt}}$ for each $\mathbf{B} \in \mathcal{B}(\mathbf{T}^{\mathrm{wt}})$, so that $\kappa$ is indeed a kernel, it suffices to observe for each bounded continuous function $F \colon \mathbf{T}^{\mathrm{wt}} \to \mathbb{R}$ that

$$\int F(\Theta(T,u,v))\mu^T(du)\nu_T(dv)$$
$$= \lim_{\varepsilon \downarrow 0} \int F(\Theta(T,u,v))\mu^{R_\varepsilon(T)}(du)\nu_T(dv)$$

and that

$$(T,d_T,\nu_T) \mapsto \int F(\Theta(T,u,v))\mu^{R_\varepsilon(T)}(du)\nu_T(dv)$$

is continuous for all $\varepsilon > 0$ (the latter follows from an argument similar to that in Lemma 7.3 of [20], where it is shown that the $(T,d_T,\nu_T) \mapsto \mu^{R_\varepsilon(T)}(T)$ is continuous). We have only sketched the argument that $\kappa$ is a kernel, because $\kappa$ is just a device for defining the measure $J$ on $\mathbf{T}^{\mathrm{wt}} \times \mathbf{T}^{\mathrm{wt}}$ in the next paragraph. It is actually the measure $J$ that we use to define our Dirichlet form, and the measure $J$ can be constructed directly as the push-forward of a measure on $U^1 \times U^1$—see the proof of Lemma 6.2.

We show in part (i) of Lemma 6.2 below that the kernel $\kappa$ is reversible with respect to the probability measure $\mathbf{P}$. More precisely, we show that if we define a measure $J$ on $\mathbf{T}^{\mathrm{wt}} \times \mathbf{T}^{\mathrm{wt}}$ by

(6.3) $$J(\mathbf{A} \times \mathbf{B}) := \int_{\mathbf{A}} \mathbf{P}(dT)\kappa(T,\mathbf{B})$$

for $\mathbf{A},\mathbf{B} \in \mathcal{B}(\mathbf{T}^{\mathrm{wt}})$, then $J$ is symmetric.

LEMMA 6.2. (i) *The measure $J$ is symmetric.*

(ii) *For each compact subset $\mathbf{K} \subset \mathbf{T}^{\mathrm{wt}}$ and open subset $\mathbf{U}$ such that $\mathbf{K} \subset \mathbf{U} \subseteq \mathbf{T}^{\mathrm{wt}}$,*

$$J(\mathbf{K}, \mathbf{T}^{\mathrm{wt}} \setminus \mathbf{U}) < \infty.$$

(iii) *The function $\Delta_{\mathrm{GH}^{\mathrm{wt}}}$ is square-integrable with respect to $J$, that is,*

$$\int_{\mathbf{T}^{\mathrm{wt}} \times \mathbf{T}^{\mathrm{wt}}} J(dT,dS)\Delta^2_{\mathrm{GH}^{\mathrm{wt}}}(T,S) < \infty.$$

PROOF. (i) Given $e',e'' \in U^1$, $0 \le u \le 1$ and $0 < \rho \le 1$, define $e^\circ(\cdot; e', e'', u, \rho) \in U^1$ by

(6.4)
$$e^\circ(t; e', e'', u, \rho)$$
$$:= \begin{cases} \mathcal{S}_{1-\rho}e''(t), & 0 \le t \le (1-\rho)u, \\ \mathcal{S}_{1-\rho}e''((1-\rho)u) + \mathcal{S}_\rho e'(t - (1-\rho)u), & \\ & (1-\rho)u \le t \le (1-\rho)u + \rho, \\ \mathcal{S}_{1-\rho}e''(t-\rho), & (1-\rho)u + \rho \le t \le 1. \end{cases}$$

That is, $e^\circ(\cdot; e', e'', u, \rho)$ is the excursion that arises from Brownian rescaling $e'$ and $e''$ to have lengths $\rho$ and $1 - \rho$, respectively, and then inserting the rescaled version of $e'$ into the rescaled version of $e''$ at a position that is a fraction $u$ of the total length of the rescaled version of $e''$.

Define a measure $\mathbb{J}$ on $U^1 \times U^1$ by

$$
\begin{aligned}
(6.5) \quad & \int_{U^1 \times U^1} \mathbb{J}(de^*, de^{**}) K(e^*, e^{**}) \\
& := \int_{[0,1]^2} du \otimes dv \, \frac{1}{2\sqrt{2\pi}} \int_0^1 \frac{d\rho}{\sqrt{(1-\rho)\rho^3}} \int \mathbb{P}(de') \otimes \mathbb{P}(de'') \\
& \quad \times K(e^\circ(\cdot; e', e'', u, \rho), e^\circ(\cdot; e', e'', v, \rho)).
\end{aligned}
$$

Clearly, the measure $\mathbb{J}$ is symmetric. It follows from the discussion at the beginning of the proof of part (i) of Lemma 5.4 and Corollary 5.2 that the measure $J$ is the push-forward of the symmetric measure $2\mathbb{J}$ by the map

$$U^1 \times U^1 \ni (e^*, e^{**}) \mapsto ((T_{2e^*}, d_{T_{2e^*}}, \nu_{T_{2e^*}}), (T_{2e^{**}}, d_{T_{2e^{**}}}, \nu_{T_{2e^{**}}})) \in \mathbf{T}^{\mathrm{wt}} \times \mathbf{T}^{\mathrm{wt}},$$

and hence $J$ is also symmetric.

(ii) The result is trivial if $\mathbf{K} = \varnothing$, so we assume that $\mathbf{K} \neq \varnothing$. Since $\mathbf{T}^{\mathrm{wt}} \setminus \mathbf{U}$ and $\mathbf{K}$ are disjoint closed sets and $\mathbf{K}$ is compact, we have that

$$(6.6) \qquad c := \inf_{T \in \mathbf{K}, S \in \mathbf{U}} \Delta_{\mathrm{GH}^{\mathrm{wt}}}(T, S) > 0.$$

Fix $T \in \mathbf{K}$. If $(u, v) \in T \times T$ is such that $\Delta_{\mathrm{GH}^{\mathrm{wt}}}(T, \Theta(T, u, v)) > c$, then $\mathrm{diam}(T) > c$ [so that we can think of $R_c(T)$, recall (2.27), as a subset of $T$]. Moreover, we claim that either:

(a) $u \in R_c(T, v)$ [recall (2.26)], or
(b) $u \notin R_c(T, v)$ and $\nu_T(S^{T,u,v}) > c$ [recall (5.17)].

Suppose, to the contrary, that $u \notin R_c(T, v)$ and that $\nu_T(S^{T,u,v}) \leq c$. Because $u \notin R_c(T, v)$, the map $f: T \to \Theta(T, u, v)$ given by

$$f(w) := \begin{cases} u, & \text{if } w \in S^{T,u,v}, \\ w, & \text{otherwise}, \end{cases}$$

is a measurable $c$-isometry. There is an analogous measurable $c$-isometry $g: \Theta(T, u, v) \to T$. Clearly,

$$d_P(f_* \nu^T, \nu^{\Theta(T,u,v)}) \leq c$$

and

$$d_P(\nu^T, g_* \nu^{\Theta(T,u,v)}) \leq c.$$

Hence, by definition, $\Delta_{\mathrm{GH}^{\mathrm{wt}}}(T, \Theta(T, u, v)) \leq c$.

Thus we have

$$J(\mathbf{K}, \mathbf{T}^{\mathrm{wt}} \setminus \mathbf{U})$$

$$\leq \int_{\mathbf{K}} \mathbf{P}\{dT\} \kappa(T, \{S : \Delta_{\mathrm{GH}^{\mathrm{wt}}}(T, S) > c\})$$

(6.7)
$$\leq \int_{\mathbf{K}} \mathbf{P}(dT) \int_T \nu_T(dv) \mu^T(R_c(T, v))$$

$$+ \int_{\mathbf{K}} \mathbf{P}(dT) \int_T \nu_T(dv) \mu^T \{u \in T : \nu_T(S^{T,u,v}) > c\}$$

$$< \infty,$$

where we have used Lemma 5.4.

(iii) Similar reasoning yields that

$$\int_{\mathbf{T}^{\mathrm{wt}} \times \mathbf{T}^{\mathrm{wt}}} J(dT, dS) \Delta_{\mathrm{GH}^{\mathrm{wt}}}^2(T, S)$$

$$= \int_{\mathbf{T}^{\mathrm{wt}}} \mathbf{P}\{dT\} \int_0^\infty dt\, 2t \kappa(T, \{S : \Delta_{\mathrm{GH}^{\mathrm{wt}}}(T, S) > t\})$$

$$\leq \int_{\mathbf{T}^{\mathrm{wt}}} \mathbf{P}(dT) \int_0^\infty dt\, 2t \int_T \nu_T(dv) \mu^T(R_c(T, v))$$

(6.8)
$$+ \int_{\mathbf{T}^{\mathrm{wt}}} \mathbf{P}(dT) \int_0^\infty dt\, 2t \int_T \nu_T(dv) \mu^T \{u \in T : \nu_T\{S^{T,u,v}\} > t\}$$

$$\leq \int_0^\infty dt\, 2t \int_{\mathbf{T}^{\mathrm{wt}}} \mathbf{P}(dT) \int_T \nu_T(dv) \mu^T(R_c(T, v))$$

$$+ \int_{\mathbf{T}^{\mathrm{wt}}} \mathbf{P}(dT) \int_T \nu_T(dv) \int_T \mu^T(du) \nu_T^2(S^{T,u,v})$$

$$< \infty,$$

where we have applied Lemma 5.4 once more. □

**7. Dirichlet forms.** Consider the bilinear form

(7.1)
$$\mathcal{E}(f, g)$$
$$:= \int_{\mathbf{T}^{\mathrm{wt}} \times \mathbf{T}^{\mathrm{wt}}} J(dT, dS)(f(S) - f(T))(g(S) - g(T)),$$

for $f, g$ in the domain

(7.2) $\mathcal{D}^*(\mathcal{E}) := \{f \in L^2(\mathbf{T}^{\mathrm{wt}}, \mathbf{P}) : f \text{ is measurable, and } \mathcal{E}(f, f) < \infty\}$

[here as usual, $L^2(\mathbf{T}^{\mathrm{wt}}, \mathbf{P})$ is equipped with the inner product $(f, g)_{\mathbf{P}} := \int \mathbf{P}(dx) \times f(x)g(x)$]. By the argument in Example 1.2.1 in [22] and Lemma 6.2, $(\mathcal{E}, \mathcal{D}^*(\mathcal{E}))$ is well defined, symmetric and Markovian.

LEMMA 7.1. *The form $(\mathcal{E}, \mathcal{D}^*(\mathcal{E}))$ is closed. That is, if $(f_n)_{n \in \mathbb{N}}$ is a sequence in $\mathcal{D}^*(\mathcal{E})$ such that*

$$\lim_{m,n \to \infty} (\mathcal{E}(f_n - f_m, f_n - f_m) + (f_n - f_m, f_n - f_m)_\mathbf{P}) = 0,$$

*then there exists $f \in \mathcal{D}^*(\mathcal{E})$ such that*

$$\lim_{n \to \infty} (\mathcal{E}(f_n - f, f_n - f) + (f_n - f, f_n - f)_\mathbf{P}) = 0.$$

PROOF. Let $(f_n)_{n \in \mathbb{N}}$ be a sequence such that $\lim_{m,n \to \infty} \mathcal{E}(f_n - f_m, f_n - f_m) + (f_n - f_m, f_n - f_m)_\mathbf{P} = 0$ [i.e., $(f_n)_{n \in \mathbb{N}}$ is Cauchy with respect to $\mathcal{E}(\cdot, \cdot) + (\cdot, \cdot)_\mathbf{P}$]. There exists a subsequence $(n_k)_{k \in \mathbb{N}}$ and $f \in L_2(\mathbf{T}^{\mathrm{wt}}, \mathbf{P})$ such that $\lim_{k \to \infty} f_{n_k} = f$, $\mathbf{P}$-a.s., and $\lim_{k \to \infty} (f_{n_k} - f, f_{n_k} - f)_\mathbf{P} = 0$. By Fatou's lemma,

$$(7.3) \qquad \int J(dT, dS)(f(S) - f(T))^2 \leq \liminf_{k \to \infty} \mathcal{E}(f_{n_k}, f_{n_k}) < \infty,$$

and so $f \in \mathcal{D}^*(\mathcal{E})$. Similarly,

$$(7.4) \qquad \begin{aligned} &\mathcal{E}(f_n - f, f_n - f) \\ &= \int J(dT, dS) \lim_{k \to \infty} ((f_n - f_{n_k})(S) - (f_n - f_{n_k})(T))^2 \\ &\leq \liminf_{k \to \infty} \mathcal{E}(f_n - f_{n_k}, f_n - f_{n_k}) \to 0 \end{aligned}$$

as $n \to \infty$. Thus $(f_n)_{n \in \mathbb{N}}$ has a subsequence that converges to $f$ with respect to $\mathcal{E}(\cdot, \cdot) + (\cdot, \cdot)_\mathbf{P}$, but, by the Cauchy property, this implies that $(f_n)_{n \in \mathbb{N}}$ itself converges to $f$. □

Let $\mathcal{L}$ denote the collection of functions $f : \mathbf{T}^{\mathrm{wt}} \to \mathbb{R}$ such that

$$(7.5) \qquad \sup_{T \in \mathbf{T}^{\mathrm{wt}}} |f(T)| < \infty$$

and

$$(7.6) \qquad \sup_{S,T \in \mathbf{T}^{\mathrm{wt}}, S \neq T} \frac{|f(S) - f(T)|}{\Delta_{\mathrm{GH}^{\mathrm{wt}}}(S, T)} < \infty.$$

Note that $\mathcal{L}$ consists of continuous functions and contains the constants. It follows from (2.16) that $\mathcal{L}$ is both a vector lattice and an algebra. By Lemma 7.2 below, $\mathcal{L} \subseteq \mathcal{D}^*(\mathcal{E})$. Therefore, the closure of $(\mathcal{E}, \mathcal{L})$ is a Dirichlet form that we will denote by $(\mathcal{E}, \mathcal{D}(\mathcal{E}))$.

LEMMA 7.2. *Suppose that $\{f_n\}_{n \in \mathbb{N}}$ is a sequence of functions from $\mathbf{T}^{\mathrm{wt}}$ into $\mathbb{R}$ such that*

$$\sup_{n \in \mathbb{N}} \sup_{T \in \mathbf{T}^{\mathrm{wt}}} |f_n(T)| < \infty,$$

$$\sup_{n \in \mathbb{N}} \sup_{S,T \in \mathbf{T}^{\mathrm{wt}}, S \neq T} \frac{|f_n(S) - f_n(T)|}{\Delta_{\mathrm{GH}^{\mathrm{wt}}}(S, T)} < \infty$$

*and*

$$\lim_{n \to \infty} f_n = f, \qquad \mathbf{P}\text{-}a.s.,$$

*for some $f : \mathbf{T}^{\mathrm{wt}} \to \mathbb{R}$. Then $\{f_n\}_{n \in \mathbb{N}} \subset \mathcal{D}^*(\mathcal{E})$, $f \in \mathcal{D}^*(\mathcal{E})$, and*

$$\lim_{n \to \infty} (\mathcal{E}(f_n - f, f_n - f) + (f_n - f, f_n - f)_{\mathbf{P}}) = 0.$$

PROOF. By the definition of the measure $J$ [see (6.3)] and the symmetry of $J$ [Lemma 6.2(i)], we have that $f_n(x) - f_n(y) \to f(x) - f(y)$ for $J$-almost every pair $(x, y)$. The result then follows from part (iii) of Lemma 6.2 and the dominated convergence theorem. □

Before showing that $(\mathcal{E}, \mathcal{D}(\mathcal{E}))$ is the Dirichlet form of a nice Markov process, we remark that $\mathcal{L}$, and hence $\mathcal{D}(\mathcal{E})$, is quite a rich class of functions: we show in the proof of Theorem 7.3 below that $\mathcal{L}$ separates points of $\mathbf{T}^{\mathrm{wt}}$ and hence if $\mathbf{K}$ is any compact subset of $\mathbf{T}^{\mathrm{wt}}$, then, by the Stone–Weierstrass theorem, the set of restrictions of functions in $\mathcal{L}$ to $\mathbf{K}$ is uniformly dense in the space of real-valued continuous functions on $\mathbf{K}$.

The following theorem states that there is a well-defined Markov process with the dynamics we would expect for a limit of the subtree prune and regraft chains.

THEOREM 7.3. *There exists a recurrent $\mathbf{P}$-symmetric Hunt process $X = (X_t, \mathbb{P}^T)$ on $\mathbf{T}^{\mathrm{wt}}$ whose Dirichlet form is $(\mathcal{E}, \mathcal{D}(\mathcal{E}))$.*

PROOF. We will check the conditions of Theorem 7.3.1 in [22] to establish the existence of $X$.

Because $\mathbf{T}^{\mathrm{wt}}$ is complete and separable (recall Theorem 2.5) there is a sequence $\mathbf{H}_1 \subseteq \mathbf{H}_2 \subseteq \cdots$ of compact subsets of $\mathbf{T}^{\mathrm{wt}}$ such that $\mathbf{P}(\bigcup_{k \in \mathbb{N}} \mathbf{H}_k) = 1$. Given $\alpha, \beta > 0$, write $\mathcal{L}_{\alpha,\beta}$ for the subset of $\mathcal{L}$ consisting of functions $f$ such that

(7.7) $$\sup_{T \in \mathbf{T}^{\mathrm{wt}}} |f(T)| \leq \alpha$$

and

(7.8) $$\sup_{S, T \in \mathbf{T}^{\mathrm{wt}}, S \neq T} \frac{|f(S) - f(T)|}{\Delta_{\mathrm{GH}^{\mathrm{wt}}}(S, T)} \leq \beta.$$

By the separability of the continuous real-valued functions on each $\mathbf{H}_k$ with respect to the supremum norm, it follows that for each $k \in \mathbb{N}$ there is a countable set $L_{\alpha,\beta,k} \subseteq \mathcal{L}_{\alpha,\beta}$ such that for every $f \in \mathcal{L}_{\alpha,\beta}$

(7.9) $$\inf_{g \in L_{\alpha,\beta,k}} \sup_{T \in \mathbf{H}_k} |f(T) - g(T)| = 0.$$

Set $L_{\alpha,\beta} := \bigcup_{k \in \mathbb{N}} L_{\alpha,\beta,k}$. Then for any $f \in \mathcal{L}_{\alpha,\beta}$ there exists a sequence $\{f_n\}_{n \in \mathbb{N}}$ in $L_{\alpha,\beta}$ such that $\lim_{n \to \infty} f_n = f$ pointwise on $\bigcup_{k \in \mathbb{N}} \mathbf{H}_k$, and hence $\mathbf{P}$-almost surely. By Lemma 7.2, the countable set $\bigcup_{m \in \mathbb{N}} L_{m,m}$ is dense in $\mathcal{L}$, and hence also dense in $\mathcal{D}(\mathcal{E})$, with respect to $\mathcal{E}(\cdot,\cdot) + (\cdot,\cdot)_\mathbf{P}$.

Now fix a countable dense subset $\mathbf{S} \subset \mathbf{T}^{\mathrm{wt}}$. Let $M$ denote the countable set of functions of the form

$$(7.10) \qquad T \mapsto p + q(\Delta_{\mathrm{GH}^{\mathrm{wt}}}(S,T) \wedge r)$$

for some $S \in \mathbf{S}$ and $p, q, r \in \mathbb{Q}$. Note that $M \subseteq \mathcal{L}$, that $M$ separates the points of $\mathbf{T}^{\mathrm{wt}}$ and, for any $T \in \mathbf{T}^{\mathrm{wt}}$, that there is certainly a function $f \in M$ with $f(T) \neq 0$.

Consequently, if $\mathcal{C}$ is the algebra generated by the countable set $M \cup \bigcup_{m \in \mathbb{N}} L_{m,m}$, then it is certainly the case that $\mathcal{C}$ is dense in $\mathcal{D}(\mathcal{E})$ with respect to $\mathcal{E}(\cdot,\cdot) + (\cdot,\cdot)_\mathbf{P}$, that $\mathcal{C}$ separates the points of $\mathbf{T}^{\mathrm{wt}}$ and, for any $T \in \mathbf{T}^{\mathrm{wt}}$, that there is a function $f \in \mathcal{C}$ with $f(T) \neq 0$.

All that remains in verifying the conditions of Theorem 7.3.1 in [22] is to check the tightness condition that there exist compact subsets $\mathbf{K}_1 \subseteq \mathbf{K}_2 \subseteq \cdots$ of $\mathbf{T}^{\mathrm{wt}}$ such that $\lim_{n \to \infty} \mathrm{Cap}(\mathbf{T}^{\mathrm{wt}} \setminus \mathbf{K}_n) = 0$, where Cap is the capacity associated with the Dirichlet form—see Remark 7.4 below for a definition. This convergence, however, is the content of Lemma 7.7 below.

Finally, because constants belong to $\mathcal{D}(\mathcal{E})$, it follows from Theorem 1.6.3 in [22] that $X$ is recurrent. $\square$

REMARK 7.4. In the proof of Theorem 7.3 we used the capacity associated with the Dirichlet form $(\mathcal{E}, \mathcal{D}(\mathcal{E}))$. We remind the reader that for an open subset $\mathbf{U} \subseteq \mathbf{T}^{\mathrm{wt}}$,

$$\mathrm{Cap}(\mathbf{U}) := \inf\{\mathcal{E}(f,f) + (f,f)_\mathbf{P} : f \in \mathcal{D}(\mathcal{E}), f(T) \geq 1, \mathbf{P}\text{-a.e. } T \in \mathbf{U}\},$$

and for a general subset $\mathbf{A} \subseteq \mathbf{T}^{\mathrm{wt}}$

$$\mathrm{Cap}(\mathbf{A}) := \inf\{\mathrm{Cap}(\mathbf{U}) : \mathbf{A} \subseteq \mathbf{U} \text{ is open}\}.$$

We refer the reader to Section 2.1 of [22] for details and a proof that Cap is a Choquet capacity.

The following results were needed in the proof of Theorem 7.3.

LEMMA 7.5. For $\varepsilon, a, \delta > 0$, put $\mathbf{V}_{\varepsilon,a} := \{T \in \mathbf{T} : \mu^T(R_\varepsilon(T)) > a\}$ and, as usual, $\mathbf{V}^\delta_{\varepsilon,a} := \{T \in \mathbf{T} : d_{\mathrm{GH}}(T, \mathbf{V}_{\varepsilon,a}) < \delta\}$. Then, for fixed $\varepsilon > 3\delta$,

$$\bigcap_{a > 0} \mathbf{V}^\delta_{\varepsilon,a} = \varnothing.$$

PROOF. Fix $S \in \mathbf{T}$. If $S \in \mathbf{V}_{\varepsilon,a}^{\delta}$, then there exists $T \in \mathbf{V}_{\varepsilon,a}$ such that $d_{\mathrm{GH}}(S,T) < \delta$. Observe that $R_{\varepsilon}(T)$ is not the trivial tree consisting of a single point because it has total length greater than $a$. Write $\{y_1, \ldots, y_n\}$ for the leaves of $R_{\varepsilon}(T)$. For all $i = 1, \ldots, n$, the connected component of $T \setminus R_{\varepsilon}(T)^o$ that contains $y_i$ contains a point $z_i$ such that $d_T(y_i, z_i) = \varepsilon$.

Let $\Re$ be a correspondence between $S$ and $T$ with $\mathrm{dis}(\Re) < 2\delta$. Pick $x_1, \ldots, x_n \in S$ such that $(x_i, z_i) \in \Re$, and hence $|d_S(x_i, x_j) - d_T(z_i, z_j)| < 2\delta$ for all $i, j$.

The distance in $R_{\varepsilon}(T)$ from the point $y_k$ to the arc $[y_i, y_j]$ is

$$(7.11) \qquad \tfrac{1}{2}(d_S(y_k, y_i) + d_S(y_k, y_j) - d_S(y_i, y_j)).$$

Thus the distance from $y_k$, $3 \leq k \leq n$, to the subtree spanned by $y_1, \ldots, y_{k-1}$ is

$$(7.12) \qquad \bigwedge_{1 \leq i \leq j \leq k-1} \tfrac{1}{2}(d_T(y_k, y_i) + d_T(y_k, y_j) - d_T(y_i, y_j)),$$

and hence

$$(7.13) \quad \begin{aligned} \mu^T(R_{\varepsilon}(T)) &= d_T(y_1, y_2) \\ &\quad + \sum_{k=3}^{n} \bigwedge_{1 \leq i \leq j \leq k-1} \tfrac{1}{2}(d_T(y_k, y_i) + d_T(y_k, y_j) - d_T(y_i, y_j)). \end{aligned}$$

Now the distance in $S$ from the point $x_k$ to the arc $[x_i, x_j]$ is

$$(7.14) \quad \begin{aligned} &\tfrac{1}{2}(d_S(x_k, x_i) + d_S(x_k, x_j) - d_S(x_i, x_j)) \\ &\geq \tfrac{1}{2}(d_T(z_k, z_i) + d_T(z_k, z_j) - d_T(z_i, z_j) - 3 \times 2\delta) \\ &= \tfrac{1}{2}(d_T(y_k, y_i) + 2\varepsilon + d_T(y_k, y_j) + 2\varepsilon - d_T(y_i, y_j) - 2\varepsilon - 6\delta) \\ &> 0 \end{aligned}$$

by the assumption that $\varepsilon > 3\delta$. In particular, $x_1, \ldots, x_n$ are leaves of the subtree spanned by $\{x_1, \ldots, x_n\}$, and $R_{\gamma}(S)$ has at least $n$ leaves when $0 < \gamma < 2\varepsilon - 6\delta$. Fix such a $\gamma$.

Now

$$(7.15) \quad \begin{aligned} &\mu^S(R_{\gamma}(S)) \\ &\geq d_S(x_1, x_2) - 2\gamma \\ &\quad + \sum_{k=3}^{n} \bigwedge_{1 \leq i \leq j \leq k-1} [\tfrac{1}{2}(d_S(x_k, x_i) + d_S(x_k, x_j) - d_S(x_i, x_j)) - \gamma] \\ &\geq \mu^T(R_{\varepsilon}(T)) + (2\varepsilon - 2\delta - 2\gamma) + (n-2)(\varepsilon - 3\delta - \gamma) \\ &\geq a + (2\varepsilon - 2\delta - 2\gamma) + (n-2)(\varepsilon - 3\delta - \gamma). \end{aligned}$$

Because $\mu^S(R_\gamma(S))$ is finite, it is apparent that $S$ cannot belong to $\mathbf{V}_{\varepsilon,a}^\delta$ when $a$ is sufficiently large. □

LEMMA 7.6. *For $\varepsilon, a > 0$, let $\mathbf{V}_{\varepsilon,a}$ be as in Lemma 7.5. Set $\mathbf{U}_{\varepsilon,a} := \{(T,\nu) \in \mathbf{T}^{\mathrm{wt}} : T \in \mathbf{V}_{\varepsilon,a}\}$. Then, for fixed $\varepsilon$,*

$$\lim_{a \to \infty} \mathrm{Cap}(\mathbf{U}_{\varepsilon,a}) = 0. \tag{7.16}$$

PROOF. Observe that $(T, d_T, \nu_T) \mapsto \mu^{R_\varepsilon(T)}(T)$ is continuous (this is essentially Lemma 7.3 of [20]), and so $\mathbf{U}_{\varepsilon,a}$ is open.

Choose $\delta > 0$ such that $\varepsilon > 3\delta$. Suppressing the dependence on $\varepsilon$ and $\delta$, define $u_a : \mathbf{T}^{\mathrm{wt}} \to [0,1]$ by

$$u_a((T,\nu)) := \delta^{-1}(\delta - d_{\mathrm{GH}}(T, \mathbf{V}_{\varepsilon,a}))_+. \tag{7.17}$$

Note that $u_a$ takes the value 1 on the open set $\mathbf{U}_{\varepsilon,a}$, and so $\mathrm{Cap}(\mathbf{U}_{\varepsilon,a}) \leq \mathcal{E}(u_a, u_a) + (u_a, u_a)_{\mathbf{P}}$. Also observe that

$$\begin{aligned}|u_a((T',\nu')) - u_a((T'',\nu''))| &\leq \delta^{-1} d_{\mathrm{GH}}(T', T'') \\ &\leq \delta^{-1} \Delta_{\mathrm{GH}^{\mathrm{wt}}}((T',\nu'),(T'',\nu'')).\end{aligned} \tag{7.18}$$

It therefore suffices by part (iii) of Lemma 6.2 and the dominated convergence theorem to show for each pair $((T',\nu'),(T'',\nu'')) \in \mathbf{T}^{\mathrm{wt}} \times \mathbf{T}^{\mathrm{wt}}$ that $u_a((T',\nu')) - u_a((T'',\nu''))$ is 0 for $a$ sufficiently large and for each $T \in \mathbf{T}^{\mathrm{wt}}$ that $u_a((T,\nu))$ is 0 for $a$ sufficiently large. However, $u_a((T',\nu')) - u_a((T'',\nu'')) \neq 0$ implies that either $T'$ or $T''$ belongs to $\mathbf{V}_{\varepsilon,a}^\delta$, while $u_a((T,\nu)) \neq 0$ implies that $T$ belongs to $\mathbf{V}_{\varepsilon,a}^\delta$. The result then follows from Lemma 7.5. □

LEMMA 7.7. *There is a sequence of compact sets $\mathbf{K}_1 \subseteq \mathbf{K}_2 \subseteq \cdots$ such that $\lim_{n \to \infty} \mathrm{Cap}(\mathbf{T}^{\mathrm{wt}} \setminus \mathbf{K}_n) = 0$.*

PROOF. By Lemma 7.6, for $n = 1, 2, \ldots$ we can choose $a_n$ so that $\mathrm{Cap}(\mathbf{U}_{2^{-n}, a_n}) \leq 2^{-n}$. Set

$$\mathbf{F}_n := \mathbf{T}^{\mathrm{wt}} \setminus \mathbf{U}_{2^{-n}, a_n} = \{(T,\nu) \in \mathbf{T}^{\mathrm{wt}} : \mu^T(R_{2^{-n}}(T)) \leq a_n\} \tag{7.19}$$

and

$$\mathbf{K}_n := \bigcap_{m \geq n} \mathbf{F}_m. \tag{7.20}$$

By Proposition 2.4 and Lemma 2.6, each set $\mathbf{K}_n$ is compact. By construction,

$$\mathrm{Cap}(\mathbf{T}^{\mathrm{wt}} \setminus \mathbf{K}_n) = \mathrm{Cap}\left(\bigcup_{m \geq n} \mathbf{U}_{2^{-m}, a_m}\right)$$

$$
\begin{aligned}
(7.21) \quad &\leq \sum_{m \geq n} \operatorname{Cap}(\mathbf{U}_{2^{-m}, a_m}) \\
&\leq \sum_{m \geq n} 2^{-m} = 2^{-(n-1)}.
\end{aligned}
$$

$\square$

**8. The trivial tree is essentially polar.** From our informal picture of the process $X$ evolving via rearrangements of the initial tree that preserve the total branch length, one might expect that if $X$ does not start at the trivial tree $T_0$ consisting of a single point, then $X$ will never hit $T_0$. However, an SPR move can decrease the diameter of a tree, so it is conceivable that, in passing to the limit, there is some probability that an infinite sequence of SPR moves will conspire to collapse the evolving tree down to a single point. Of course, it is hard to imagine from the approximating dynamics how $X$ could recover from such a catastrophe—which it would have to since it is reversible with respect to the continuum random tree distribution.

In this section we will use potential theory for Dirichlet forms to show that $X$ does not hit $T_0$ from $\mathbf{P}$-almost all starting points; that is, that the set $\{T_0\}$ is *essentially polar*.

Let $\bar{d}$ be the map which sends a weighted $\mathbb{R}$ tree $(T, d, \nu)$ to the $\nu$-averaged distance between pairs of points in $T$. That is,

$$
(8.1) \quad \bar{d}((T, d, \nu)) := \int_T \int_T \nu(dx)\nu(dy)\, d(x, y), \qquad (T, d, \nu) \in \mathbf{T}^{\mathrm{wt}}.
$$

In order to show that $T_0$ is essentially polar, it will suffice to show that the set

$$
(8.2) \quad \{(T, d, \nu) \in \mathbf{T}^{\mathrm{wt}} : \bar{d}((T, d, \nu)) = 0\}
$$

is essentially polar.

LEMMA 8.1. *The function $\bar{d}$ belongs to the domain $\mathcal{D}(\mathcal{E})$.*

PROOF. If we let $\bar{d}_n((T, d, \nu)) := \int_T \int_T \nu(dx)\nu(dy)[d(x, y) \wedge n]$, for $n \in \mathbb{N}$, then $\bar{d}_n \uparrow \bar{d}$, $\mathbf{P}$-a.s. By the triangle inequality,

$$
(8.3) \quad (\bar{d}, \bar{d})_\mathbf{P} \leq \int \mathbf{P}(dT)(\operatorname{diam}(T))^2 \leq \int \mathbb{P}(de)\left(4 \sup_{t \in [0,1]} e(t)\right)^2 < \infty,
$$

and hence $\bar{d}_n \to \bar{d}$ as $n \to \infty$ in $L^2(\mathbf{T}^{\mathrm{wt}}, \mathbf{P})$.

Notice, moreover, that for $(T, d, \nu) \in \mathbf{T}^{\mathrm{wt}}$ and $u, v \in T$,

$$
\begin{aligned}
&(\bar{d}((T, d, \nu)) - \bar{d}(\Theta((T, d, \nu), u, v)))^2 \\
(8.4) \quad &= 2 \int_{S^{T,u,v}} \int_{T \setminus S^{T,u,v}} \nu(dx)\nu(dy)(d(y, u) - d(y, v))^2 \\
&= 2\nu_T(S^{T,u,v})\nu(T \setminus S^{T,u,v})\, d^2(u, v).
\end{aligned}
$$

Hence, applying Corollary 5.2 and the invariance of the standard Brownian excursion under random rerooting (see Section 2.7 of [3]),

$$\int_{\mathbf{T}^{\mathrm{wt}} \times \mathbf{T}^{\mathrm{wt}}} J(dT, dS)(\bar{d}(T) - \bar{d}(S))^2$$

$$= 2 \int_{\mathbf{T}^{\mathrm{wt}}} \mathbf{P}(dT)$$

$$\times \int_{T \times T} \nu_T(dv) \mu^T(du) \nu_T(S^{T,u,v}) \nu_T(T \setminus S^{T,u,v}) d_T^2(u,v)$$

(8.5) $$\leq 2 \int \mathbb{P}(de) 2 \int_{\Gamma_e} \frac{ds \otimes da}{\bar{s}(e,s,a) - \underline{s}(e,s,a)} \zeta(\hat{e}^{s,a}) \zeta(\check{e}^{s,a})(2a)^2$$

$$= \frac{8}{\sqrt{2\pi}} \int_0^1 \frac{d\rho}{\sqrt{(1-\rho)\rho^3}}$$

$$\times \int \mathbb{P}(de') \otimes \mathbb{P}(de'') \rho(1-\rho)(\sup \mathcal{S}_{1-\rho} e'')^2$$

$$= \frac{8}{\sqrt{2\pi}} \int_0^1 \frac{d\rho}{\sqrt{(1-\rho)\rho^3}} \rho(1-\rho)^2 \int \mathbb{P}(de) \left( \sup_{t \in [0,1]} e(t) \right)^2 < \infty.$$

Consequently, by dominated convergence, $\mathcal{E}(\bar{d} - \bar{d}_n, \bar{d} - \bar{d}_n) \to 0$ as $n \to \infty$.

It is therefore enough to verify that $\bar{d}_n \in \mathcal{L}$ for all $n \in \mathbb{N}$. Obviously,

(8.6) $$\sup_{T \in \mathbf{T}^{\mathrm{wt}}} \bar{d}_n(T) \leq n,$$

and so the boundedness condition (7.5) holds. To show that the "Lipschitz" property (7.6) holds, fix $\varepsilon > 0$, and let $(T, \nu_T), (S, \nu_S) \in \mathbf{T}^{\mathrm{wt}}$ be such that $\Delta_{\mathrm{GH}^{\mathrm{wt}}}((T, \nu_T), (S, \nu_S)) < \varepsilon$. Then there exist $f \in F_{T,S}^\varepsilon$ and $g \in F_{S,T}^\varepsilon$ such that $d_\mathrm{P}(\nu_T, g_* \nu_S) < \varepsilon$ and $d_\mathrm{P}(f_* \nu_T, \nu_S) < \varepsilon$ [recall $F_{T,S}^\varepsilon$ from (2.10)]. Hence

$$|\bar{d}_n((T, \nu_T)) - \bar{d}_n((S, \nu_S))|$$

$$\leq \left| \int_T \int_T \nu_T(dx) \nu_T(dy) (d_T(x,y) \wedge n) \right.$$

(8.7) $$\left. - \int_{g(S)} \int_{g(S)} g_* \nu_S(dx) g_* \nu_S(dy) (d_T(x,y) \wedge n) \right|$$

$$+ \left| \int_{g(S)} \int_{g(S)} g_* \nu_S(dx) g_* \nu_S(dy) (d_T(x,y) \wedge n) \right.$$

$$\left. - \int_S \int_S \nu_S(dx') \nu_S(dy') (d_S(x',y') \wedge n) \right|.$$

For the first term on the right-hand side of (8.7) we get

$$\left| \int_T \int_T \nu_T(dx) \nu_T(dy) (d_T(x,y) \wedge n) \right.$$

$$
\begin{aligned}
&\left. - \int_{g(S)} \int_{g(S)} g_*\nu_S(dx) g_*\nu_S(dy)(d_T(x,y) \wedge n) \right| \\
&\leq \left| \int_T \int_T \nu_T(dx)\nu_T(dy)(d_T(x,y) \wedge n) \right. \\
&\quad \left. - \int_T \int_{g(S)} \nu_T(dx) g_*\nu_S(dy)(d_T(x,y) \wedge n) \right| \\
&\quad + \left| \int_{S(g)} \int_T g_*\nu_S(dx)\nu_T(dy)(d_T(x,y) \wedge n) \right. \\
&\quad \left. - \int_{g(S)} \int_{g(S)} g_*\nu_S(dx) g_*\nu_S(dy)(d_T(x,y) \wedge n) \right|.
\end{aligned}
\tag{8.8}
$$

By assumption and Theorem 3.1.2 in [19], we can find a probability measure $\nu$ on $T \times T$ with marginals $\nu_T$ and $g_*\nu_S$ such that

$$
\nu\{(x,y) : d_T(x,y) \geq \varepsilon\} \leq \varepsilon.
\tag{8.9}
$$

Hence, for all $x \in T$,

$$
\begin{aligned}
&\left| \int_T \nu_T(dy)(d_T(x,y) \wedge n) - \int_{g(S)} g_*\nu_S(dy)(d_T(x,y) \wedge n) \right| \\
&\leq \int_{T \times g(S)} \nu(d(y,y'))|(d_T(x,y) \wedge n) - (d_T(x,y') \wedge n)| \\
&\leq \int_{T \times g(S)} \nu(d(y,y'))(d_T(y,y') \wedge n) \\
&\leq (1 + (\mathrm{diam}(T) \wedge n)) \cdot \varepsilon.
\end{aligned}
\tag{8.10}
$$

For the second term in (8.7) we use the fact that $g$ is an $\varepsilon$-isometry, that is, $|(d_S(x',y') \wedge n) - (d_T(g(x'), g(y')) \wedge n)| < \varepsilon$ for all $x', x'' \in T$. A change of variables then yields that

$$
\begin{aligned}
&\left| \int_{g(S)} \int_{g(S)} g_*\nu_S(dx) g_*\nu_S(dy)(d_T(x,y) \wedge n) \right. \\
&\quad \left. - \int_S \int_S \nu_S(dx')\nu_S(dy')(d_S(x',y') \wedge n) \right| \\
&\leq \varepsilon + \left| \int_{g(S)} \int_{g(S)} g_*\nu_S(dx) g_*\nu_S(dy)(d_T(x,y) \wedge n) \right. \\
&\quad \left. - \int_S \int_S \nu_S(dx')\nu_S(dy')(d_T(g(x'), g(y')) \wedge n) \right| \\
&= \varepsilon.
\end{aligned}
\tag{8.11}
$$

Combining (8.7) through (8.11) yields finally that

$$(8.12) \quad \sup_{(T,\nu_T) \neq (S,\nu_S) \in \mathbf{T}^{\mathrm{wt}}} \frac{|\bar{d}_n((T,\nu_T)) - \bar{d}_n((S,\nu_S))|}{\Delta_{\mathrm{GH}^{\mathrm{wt}}}((T,\nu_T),(S,\nu_S))} \leq 3 + 2n. \qquad \square$$

PROPOSITION 8.2. *The set $\{T \in \mathbf{T}^{\mathrm{wt}} : \bar{d}(T) = 0\}$ is essentially polar. In particular, the set $\{T_0\}$ consisting of the trivial tree is essentially polar.*

PROOF. We need to show that $\mathrm{Cap}(\{T \in \mathbf{T}^{\mathrm{wt}} : \bar{d}(T) = 0\}) = 0$ (see Theorem 4.2.1 of [22]).

For $\varepsilon > 0$ set

$$(8.13) \quad \mathbf{W}_\varepsilon := \{T \in \mathbf{T}^{\mathrm{wt}} : \bar{d}(T) < \varepsilon\}.$$

By the argument in the proof of Lemma 8.1, the function $\bar{d}$ is continuous, and so $\mathbf{W}_\varepsilon$ is open. It suffices to show that $\mathrm{Cap}(\mathbf{W}_\varepsilon) \downarrow 0$ as $\varepsilon \downarrow 0$.

Put

$$(8.14) \quad u_\varepsilon(T) := \left(2 - \frac{\bar{d}(T)}{\varepsilon}\right)_+, \qquad T \in \mathbf{T}^{\mathrm{wt}}.$$

Then $u \in \mathcal{D}(\mathcal{E})$ by Lemma 8.1 and the fact that the domain of a Dirichlet form is closed under composition with Lipschitz functions. Because $u_\varepsilon(T) \geq 1$ for $T \in \mathbf{W}_\varepsilon$, it thus further suffices to show

$$(8.15) \quad \lim_{\varepsilon \downarrow 0}(\mathcal{E}(u_\varepsilon, u_\varepsilon) + (u_\varepsilon, u_\varepsilon)_{\mathbf{P}}) = 0.$$

By elementary properties of the standard Brownian excursion,

$$(8.16) \quad (u_\varepsilon, u_\varepsilon)_{\mathbf{P}} \leq 4\mathbf{P}\{T : \bar{d}(T) < 2\varepsilon\} \to 0$$

as $\varepsilon \downarrow 0$. Estimating $\mathcal{E}(u_\varepsilon, u_\varepsilon)$ will be somewhat more involved.

Let $\hat{E}$ and $\check{E}$ be two independent standard Brownian excursions, and let $U$ and $V$ be two independent random variables that are independent of $\hat{E}$ and $\check{E}$ and uniformly distributed on $[0,1]$. With a slight abuse of notation, we will write $\mathbb{P}$ for the probability measure on the probability space where $\hat{E}$, $\check{E}$, $U$ and $V$ are defined.

Set

$$\hat{D} := 4\int_{0 \leq s < t \leq 1} ds \otimes dt\left[\hat{E}_s + \hat{E}_t - 2\inf_{s \leq w \leq t}\hat{E}_w\right],$$

$$\hat{H} := 2\int_{[0,1]} dt\hat{E}_t,$$

$$(8.17) \quad \check{D} := 4\int_{0 \leq s < t \leq 1} ds \otimes dt\left[\check{E}_s + \check{E}_t - 2\inf_{s \leq w \leq t}\check{E}_w\right],$$

$$\check{H}_U := 2\int_{[0,1]} dt\left[\check{E}_t + \check{E}_U - 2\inf_{U\wedge t\leq w\leq U\vee t}\check{E}_w\right],$$

$$\check{H}_V := 2\int_{[0,1]} dt\left[\check{E}_t + \check{E}_V - 2\inf_{V\wedge t\leq w\leq V\vee t}\check{E}_w\right].$$

For $0\leq \rho\leq 1$ set

(8.18)
$$\begin{aligned}D_U(\rho) &:= (1-\rho)^2\sqrt{1-\rho}\check{D} + \rho^2\sqrt{\rho}\hat{D}\\ &\quad + 2(1-\rho)\rho\sqrt{\rho}\hat{H} + 2(1-\rho)\rho\sqrt{1-\rho}\check{H}_U\end{aligned}$$

and

(8.19)
$$\begin{aligned}D_V(\rho) &:= (1-\rho)^2\sqrt{1-\rho}\check{D} + \rho^2\sqrt{\rho}\hat{D}\\ &\quad + 2(1-\rho)\rho\sqrt{\rho}\hat{H} + 2(1-\rho)\rho\sqrt{1-\rho}\check{H}_V.\end{aligned}$$

Then

(8.20)
$$\begin{aligned}&\mathcal{E}(u_\varepsilon, u_\varepsilon)\\ &= \frac{1}{2\sqrt{2\pi}}\mathbb{P}\left[\int_0^1 \frac{d\rho}{\sqrt{(1-\rho)\rho^3}}\right.\\ &\quad\left.\times\left\{\left(2-\frac{D_U(\rho)}{\varepsilon}\right)_+ - \left(2-\frac{D_V(\rho)}{\varepsilon}\right)_+\right\}^2\right].\end{aligned}$$

Fix $0 < a < \frac{1}{2}$ and write $\bar{a} = 1-a$ for convenience. We can write the right-hand side of (8.20) as the sum of three terms $I(\varepsilon, a)$, $II(\varepsilon, a)$ and $III(\varepsilon, a)$, that arise from integrating $\rho$ over the respective ranges

(8.21) $$\{\rho : D_U(\rho)\vee D_V(\rho) \leq 2\varepsilon, 0\leq \rho\leq a\},$$

(8.22) $$\{\rho : D_U(\rho)\wedge D_V(\rho) \leq 2\varepsilon \leq D_U(\rho)\vee D_V(\rho), 0\leq \rho\leq a\},$$

and

(8.23) $$\{\rho : a < \rho \leq 1\}.$$

Consider $I(\varepsilon, a)$ first. Note that if $D_U(\rho)\vee D_V(\rho) \leq 2\varepsilon$, then

(8.24) $$\left\{\left(2-\frac{D_U(\rho)}{\varepsilon}\right)_+ - \left(2-\frac{D_V(\rho)}{\varepsilon}\right)_+\right\}^2 \leq 2^2\frac{\rho^2}{\varepsilon^2}\{\check{H}_U - \check{H}_V\}^2.$$

Moreover,

(8.25)
$$\begin{aligned}&\{0\leq \rho\leq a : D_U(\rho)\vee D_V(\rho) \leq 2\varepsilon\}\\ &\quad \subseteq \{0\leq \rho\leq a : (1-\rho)^{5/2}\check{D} + 2(1-\rho)^{3/2}\rho(\check{H}_U\vee\check{H}_V) \leq 2\varepsilon\}\\ &\quad \subseteq \{0\leq \rho\leq a : \bar{a}^{5/2}\check{D} + 2\bar{a}^{3/2}\rho(\check{H}_U\vee\check{H}_V) \leq 2\varepsilon\}\\ &\quad = \left\{\rho : 0\leq \rho\leq \frac{(2\varepsilon - \bar{a}^{5/2}\check{D})_+}{2\bar{a}^{3/2}(\check{H}_U\vee\check{H}_V)}\wedge a\right\}.\end{aligned}$$

Thus $I(\varepsilon, a)$ is bounded above by the expectation of the random variable that arises from integrating $2^2\rho^2\{\check{H}_U - \check{H}_V\}^2/\varepsilon^2$ against the measure $\frac{1}{2\sqrt{2\pi}}\frac{d\rho}{\sqrt{(1-\rho)\rho^3}}$ over the interval $[0, (2\varepsilon - \bar{a}^{5/2}\check{D})_+/(2\bar{a}^{3/2}(\check{H}_U \vee \check{H}_V))]$. Note that

$$\text{(8.26)} \qquad \int_0^x \frac{d\rho}{\sqrt{\rho^3}}\rho^\alpha = \frac{1}{\alpha - 1/2}x^{\alpha-1/2}, \qquad \alpha > \frac{1}{2}.$$

Hence, letting $C$ denote a generic constant with a value that does not depend on $\varepsilon$ or $a$ and may change from line to line,

$$\text{(8.27)} \qquad \begin{aligned} I(\varepsilon, a) &\leq C\mathbb{P}\left[\left(\frac{(2\varepsilon - \bar{a}^{5/2}\check{D})_+}{\check{H}_U \vee \check{H}_V}\right)^{3/2} \frac{\{\check{H}_U - \check{H}_V\}^2}{\varepsilon^2}\right] \\ &\leq \frac{C}{\varepsilon^2}\mathbb{P}[(2\varepsilon - \bar{a}^{5/2}\check{D})_+^{3/2}(\check{H}_U \vee \check{H}_V)^{1/2}] \\ &\leq \frac{C}{\varepsilon^{1/2}}\mathbb{P}[(\check{H}_U + \check{H}_V)^{1/2}\mathbb{1}\{\check{D} \leq 2\bar{a}^{-5/2}\varepsilon\}] \\ &\leq \frac{C}{\varepsilon^{1/2}}\mathbb{P}[\check{D}^{1/2}\mathbb{1}\{\check{D} \leq 2\bar{a}^{-5/2}\varepsilon\}] \\ &\leq C\mathbb{P}\{\check{D} \leq 2\bar{a}^{-5/2}\varepsilon\}, \end{aligned}$$

where in the second last line we used the fact that

$$\text{(8.28)} \qquad \mathbb{P}[\check{H}_U|\check{E}] = \mathbb{P}[\check{H}_V|\check{E}] = \check{D},$$

and Jensen's inequality for conditional expectations to obtain the inequalities $\mathbb{P}[\check{H}_U^{1/2}|\check{E}] \leq \check{D}^{1/2}$ and $\mathbb{P}[\check{H}_V^{1/2}|\check{E}] \leq \check{D}^{1/2}$. Thus, $\lim_{\varepsilon \downarrow 0} I(\varepsilon, a) = 0$ for any value of $a$.

Turning to $II(\varepsilon, a)$, first note that $\hat{D} \leq 4\hat{H}$ and, by the triangle inequality,

$$\text{(8.29)} \qquad \check{D} \leq 2(\check{H}_U \wedge \check{H}_V).$$

Hence, for some constant $K$ that does not depend on $\varepsilon$ or $a$,

$$\text{(8.30)} \qquad |D_U(\rho) \wedge D_V(\rho) - \check{D}| \leq K(\hat{H}\rho^{3/2} + (\check{H}_U \wedge \check{H}_V)\rho)$$

and

$$\text{(8.31)} \qquad |D_U(\rho) \vee D_V(\rho) - \check{D}| \leq K(\hat{H}\rho^{3/2} + (\check{H}_U \vee \check{H}_V)\rho).$$

Combining (8.31) with an argument similar to that which established (8.25) gives, for a suitable constant $K^*$,

$$\begin{aligned} &\{0 \leq \rho \leq a : D_U(\rho) \wedge D_V(\rho) \leq 2\varepsilon \leq D_U(\rho) \vee D_V(\rho)\} \\ &= \{0 \leq \rho \leq a : 2\varepsilon \leq D_U(\rho) \vee D_V(\rho)\} \\ &\quad \cap \{0 \leq \rho \leq a : D_U(\rho) \wedge D_V(\rho) \leq 2\varepsilon\} \end{aligned} \tag{8.32}$$

$$\subseteq \left\{\rho: \frac{(2\varepsilon - \check{D})_+}{K^*(\hat{H} + \check{H}_U \vee \check{H}_V)} \leq \rho \leq a\right\}$$

$$\cap \left\{\rho: 0 \leq \rho \leq \frac{(2\varepsilon - \bar{a}^{5/2}\check{D})_+}{2\bar{a}^{3/2}(\check{H}_U \wedge \check{H}_V)} \wedge a\right\}.$$

Moreover, by (8.30) and the observation $|(2\varepsilon - x)_+ - (2\varepsilon - y)_+| \leq |x - y|$, we have for $D_U(\rho) \wedge D_V(\rho) \leq 2\varepsilon \leq D_U(\rho) \vee D_V(\rho)$ that

$$\left\{\left(2 - \frac{D_U(\rho)}{\varepsilon}\right)_+ - \left(2 - \frac{D_V(\rho)}{\varepsilon}\right)_+\right\}^2$$

$$= \left\{\left(2 - \frac{D_U(\rho) \wedge D_V(\rho)}{\varepsilon}\right)_+\right\}^2$$

$$\leq \frac{2}{\varepsilon^2}\{(2\varepsilon - \check{D})_+\}^2$$

(8.33)

$$+ \frac{2}{\varepsilon^2}\{(2\varepsilon - D_U(\rho) \wedge D_V(\rho))_+ - (2\varepsilon - \check{D})_+\}^2$$

$$\leq \frac{2}{\varepsilon^2}\{(2\varepsilon - \check{D})_+\}^2 + \frac{2}{\varepsilon^2}\{D_U(\rho) \wedge D_V(\rho) - \check{D}\}^2$$

$$\leq \frac{C}{\varepsilon^2}[(2\varepsilon - \check{D})_+^2 + \hat{H}^2 \rho^3 + (\check{H}_U \wedge \check{H}_V)^2 \rho^2],$$

for a suitable constant $C$ that does not depend on $\varepsilon$ or $a$. It follows from (8.26) and

(8.34) $$\int_x^a \frac{d\rho}{\sqrt{\rho^3}} \rho^\beta = \frac{1}{1/2 - \beta}[x^{\beta - 1/2} - a^{\beta - 1/2}], \qquad \beta < \frac{1}{2},$$

that

(8.35)
$$II(\varepsilon, a) \leq \frac{C'}{\varepsilon^2}\mathbb{P}\left[(2\varepsilon - \check{D})_+^2\left\{\frac{(2\varepsilon - \check{D})_+}{\hat{H} + \check{H}_U \vee \check{H}_V}\right\}^{-1/2}\right]$$

$$+ \frac{C''}{\varepsilon^2}\mathbb{P}\left[\hat{H}^2\left\{\frac{(2\varepsilon - \bar{a}^{5/2}\check{D})_+}{2\bar{a}^{3/2}(\check{H}_U \wedge \check{H}_V)} \wedge a\right\}^{5/2}\right]$$

$$+ \frac{C'''}{\varepsilon^2}\mathbb{P}\left[(\check{H}_U \wedge \check{H}_V)^2\left\{\frac{(2\varepsilon - \bar{a}^{5/2}\check{D})_+}{2\bar{a}^{3/2}(\check{H}_U \wedge \check{H}_V)} \wedge a\right\}^{3/2}\right]$$

for suitable constants $C'$, $C''$ and $C'''$.

Consider the first term in (8.35). Using Jensen's inequality for conditional expectations and (8.28) again, this term is bounded above by

(8.36) $$\frac{1}{\varepsilon^2}\mathbb{P}[(2\varepsilon - \check{D})_+^{3/2}\{A\check{D}^{1/2} + B\}] \leq \frac{1}{\varepsilon^2}\mathbb{P}[(2\varepsilon - \check{D})_+^{3/2}\{2^{1/2}A\varepsilon^{1/2} + B\}]$$

for suitable constants $A, B$. Now, by Jensen's inequality for conditional expectation yet again, along with the invariance of standard Brownian excursion under random rerooting (see Section 2.7 of [3]) and the fact that

$$\mathbb{P}\{\check{E}_U \in dr\} = re^{-r^2/2} \, dr \tag{8.37}$$

(see Section 3.3 of [3]), we have

$$\mathbb{P}[(2\varepsilon - \check{D})_+^{3/2}]$$

$$= \mathbb{P}\left[\left(\mathbb{P}\left[2\varepsilon - 2\left\{\check{E}_U + \check{E}_V - 2\inf_{U \wedge V \leq t \leq U \vee V} \check{E}_t\right\} \Big| \check{E}\right]\right)_+^{3/2}\right]$$

$$\leq \mathbb{P}\left[\left(2\varepsilon - 2\left\{\check{E}_U + \check{E}_V - 2\inf_{U \wedge V \leq t \leq U \vee V} \check{E}_t\right\}\right)_+^{3/2}\right] \tag{8.38}$$

$$= \mathbb{P}[(2\varepsilon - 2\check{E}_U)_+^{3/2}] = \int_0^\infty dr \, re^{-r^2/2}(2\varepsilon - 2r)_+^{3/2}$$

$$\leq \int_0^\varepsilon dr \, r(2\varepsilon - 2r)^{3/2} = 2^{3/2}\varepsilon^{7/2} \int_0^1 ds \, s(1-s)^{3/2}.$$

Thus the limit as $\varepsilon \downarrow 0$ of the first term in (8.35) is 0 for each $a$.

For the second term in (8.35), first observe by Jensen's inequality for conditional expectation and (8.37) that

$$\mathbb{P}\{\check{D} \leq r\} \leq \mathbb{P}\left[\left(2 - \frac{\check{D}}{r}\right)_+\right] \leq \mathbb{P}\left[\left(2 - \frac{\check{E}_U}{r}\right)_+\right]$$

$$\leq 2\mathbb{P}\{\check{E}_U \leq 2r\} \leq 2\frac{(2r)^2}{2} = 4r^2. \tag{8.39}$$

Combining this observation with (8.29) and integrating by parts gives

$$\mathbb{P}\left[\left\{\frac{(2\varepsilon - \bar{a}^{5/2}\check{D})_+}{2\bar{a}^{3/2}(\check{H}_U \wedge \check{H}_V)} \wedge a\right\}^{5/2}\right]$$

$$\leq \mathbb{P}\left[\left\{\frac{(2\varepsilon - \bar{a}^{5/2}\check{D})_+}{\bar{a}^{3/2}\check{D}} \wedge a\right\}^{5/2}\right]$$

$$= \int_0^{2\varepsilon/\bar{a}^{5/2}} \mathbb{P}\{\check{D} \in dr\}\left(\bar{a}^{-3/2}\left(\frac{2\varepsilon}{r} - \bar{a}^{5/2}\right) \wedge a\right)^{5/2} \tag{8.40}$$

$$\leq \int_{2\varepsilon/(a\bar{a}^{1/2}+\bar{a}^{5/2})}^{2\varepsilon/\bar{a}^{5/2}} dr \, 4r^2 \bar{a}^{-15/4} \frac{5}{2}\left(\frac{2\varepsilon}{r} - \bar{a}^{5/2}\right)^{3/2} \frac{2\varepsilon}{r^2}$$

$$= 40\varepsilon^2 \bar{a}^{-15/4} \int_{1/(a\bar{a}^{1/2}+\bar{a}^{5/2})}^{1/\bar{a}^{5/2}} ds\left(\frac{1}{s} - \bar{a}^{5/2}\right)^{3/2}.$$

SUBTREE PRUNE AND REGRAFT 43

If we denote the rightmost term by $L(\varepsilon, a)$, then it is clear that

$$\lim_{a \downarrow 0} \lim_{\varepsilon \downarrow 0} \frac{1}{\varepsilon^2} L(\varepsilon, a) = 0. \tag{8.41}$$

From (8.28) and Jensen's inequality for conditional expectations, the third term in (8.35) is bounded above by

$$\frac{C}{\varepsilon^2} \mathbb{P}[(\check{H}_U \wedge \check{H}_V)^{1/2}(2\varepsilon - \bar{a}^{5/2}\check{D})_+^{3/2}] \leq \frac{C}{\varepsilon^2} \mathbb{P}[\check{D}^{1/2}(2\varepsilon - \bar{a}^{5/2}\check{D})_+^{3/2}] \\ \leq \frac{C}{\varepsilon^{3/2}} \mathbb{P}[(2\varepsilon - \bar{a}^{5/2}\check{D})_+^{3/2}], \tag{8.42}$$

and the calculation in (8.38) shows that the rightmost term converges to zero as $\varepsilon \downarrow 0$ for each $a$.

Putting together the observations we have made on the three terms in (8.35), we see that

$$\lim_{a \downarrow 0} \lim_{\varepsilon \downarrow 0} II(\varepsilon, a) = 0. \tag{8.43}$$

It follows from the dominated convergence theorem that

$$\lim_{\varepsilon \downarrow 0} III(\varepsilon, a) = 0 \tag{8.44}$$

for all $a$, and this completes the proof. $\square$

**Acknowledgments.** The authors thank David Aldous, Peter Pfaffelhuber, Jim Pitman and Tandy Warnow for helpful discussions, and the referee and an Associate Editor for extremely close readings of the paper and several very helpful suggestions. Part of the research was conducted while the authors were visiting the Pacific Institute for the Mathematical Sciences in Vancouver, Canada.## REFERENCES

[1] ABRAHAM, R. and SERLET, L. (2002). Poisson snake and fragmentation. *Electron. J. Probab.* **7** 15–17 (electronic). MR1943890
[2] ALDOUS, D. (1991). The continuum random tree. I. *Ann. Probab.* **19** 1–28. MR1085326
[3] ALDOUS, D. (1991). The continuum random tree. II. An overview. In *Stochastic Analysis* (M. T. Barlow and N. H. Bingham, eds.) 23–70. Cambridge Univ. Press. MR1166406
[4] ALDOUS, D. (1993). The continuum random tree III. *Ann. Probab.* **21** 248–289. MR1207226
[5] ALDOUS, D. J. (2000). Mixing time for a Markov chain on cladograms. *Combin. Probab. Comput.* **9** 191–204. MR1774749
[6] ALLEN, B. L. and STEEL, M. (2001). Subtree transfer operations and their induced metrics on evolutionary trees. *Ann. Comb.* **5** 1–15. MR1841949

DEPARTMENT OF STATISTICS  
UNIVERSITY OF CALIFORNIA AT BERKELEY  
367 EVANS HALL  
BERKELEY, CALIFORNIA 94720-3860  
USA  
E-MAIL: evans@stat.berkeley.edu  
URL: www.stat.berkeley.edu/users/evans/

MATHEMATISCHES INSTITUT  
UNIVERSITÄT ERLANGEN–NÜRNBERG  
BISMARCKSTRASSE $1\frac{1}{2}$  
91054 ERLANGEN  
GERMANY  
E-MAIL: winter@mi.uni-erlangen.de  
URL: www.mi.uni-erlangen.de/~winter/